\documentclass[10pt,twocolumn,twoside]{IEEEtran}
\usepackage{bbm}
\usepackage{amsmath}
\usepackage{amsfonts}
\usepackage{mathrsfs}
\usepackage{amssymb}
\usepackage{multirow}
\usepackage{latexsym}
\usepackage{psfrag}
\usepackage{graphicx}
\usepackage{cite}
\usepackage{color}
\usepackage{times,balance}
\usepackage{algorithm,multirow,xcolor}
\usepackage{algorithmicx}
\usepackage{algpseudocode}
\usepackage{longtable}
\usepackage{booktabs}
\usepackage{colortbl}
\usepackage{bm}
\usepackage{booktabs}
\usepackage{amssymb}
\usepackage{pifont}
\usepackage{wasysym}
\usepackage{threeparttable}

\newcommand{\cmark}{\ding{51}}%
\newcommand{\xmark}{\ding{55}}%
\definecolor{mygray}{gray}{.8}

\newtheorem{theorem}{Theorem}

\newtheorem{assumption}{Assumption}

\newtheorem{lemma}{Lemma}
\newtheorem{remark}{Remark}

\allowdisplaybreaks[4]

\begin{document}


\title{Distributed Online Optimization in Time-Varying Unbalanced Networks without Explicit
Subgradients}

%



\author{Yongyang Xiong, Xiang Li, Keyou You, \IEEEmembership{Senior Member, IEEE},
and Ligang Wu, \IEEEmembership{Fellow, IEEE}
\thanks{This work was supported in part by the National Natural Science Foundation
of China (61525303, 41772377 and 61673130), Top-Notch Young Talents Program of China
(Ligang Wu), and the Self-Planned Task of State Key Laboratory of Robotics
and System (HIT) (SKLRS201806B).
\emph{(Corresponding author: Ligang Wu.)}}
\thanks{
Y. Xiong, and K. You are with the 
Department of Automation, Beijing National Research Center for Information Science and Technology, Tsinghua University, Beijing 100084, P. R. China. This work was done when Y. Xiong was with the
Department of Control Science and Engineering, Harbin Institute of Technology. 
E-mail: \texttt{xiongyy@tsinghua.edu.cn}; \texttt{youky@tsinghua.edu.cn}.}
\thanks{X. Li and L. Wu are with the
Department of Control Science and Engineering, Harbin Institute of Technology, Harbin 150001, P. R. China.
E-mail: \texttt{lixiang\_soa@hit.edu.cn}; \texttt{ligangwu@hit.edu.cn}.}}


\maketitle

\begin{abstract}
This paper studies a distributed online constrained optimization problem over time-varying unbalanced digraphs without explicit subgradients. In sharp contrast to the existing algorithms, we design a novel consensus-based distributed online algorithm with a local  randomized zeroth-order oracle and then rescale the oracle by constructing row-stochastic matrices, which aims to address the unbalancedness of time-varying digraphs. Under mild conditions, 
the average dynamic  regret over a time horizon is shown to asymptotically converge 
at a sublinear rate 
provided that the accumulated variation grows sublinearly with a specific order. 
Moreover, the counterpart of the proposed algorithm when subgradients are available is also provided, along with its dynamic regret bound, which reflects that the convergence of our algorithm is essentially not affected by the zeroth-order oracle. Simulations on distributed targets tracking problem and dynamic sparse signal recovery problem in sensor networks are employed to demonstrate the effectiveness of the proposed algorithm.
\end{abstract}


\section{Introduction}

\IEEEPARstart{D}{istributed}
optimization has been undoubtedly attracting a surge of
attentions in recent years with the rapid development of large-scale networks.
Many practical issues can be solved
within the framework of distributed optimization problems,
such as 
privacy preserving \cite{LiQ2020Privacy,XiongY2020Privacy}, 
resource allocation \cite{ZhangJ2019Distributed,XuJ2019A,XuY2017A}, 
sensor fusion \cite{ZhuS2018Mitigating}, 
just to mention a few. 
{\color{black}These applications promote the design of distributed algorithms such that a group of nodes 
cooperatively optimize the sum of their local cost functions via local communications.} 
See \cite{YouK2019Distributed,ReisizadehA2019An,Jakovetic2014Convergence}  and the references therein.

%
Despite that distributed optimization algorithms has been studied extensively, {\color{black}the dynamic aspect of the problem has not been fully addressed where the cost functions may vary with time in an uncertain and even adversarial fashion}.
  For example, 
{\color{black}time-varying cost functions frequently appear in machine learning, 
where data samples are observed in a sequential manner and the newly observed data samples result in new cost functions.}  
   This inspires us to extend the distributed optimization algorithms to an online setting,
in which the local cost functions vary with time and are only revealed to individuals after each node
has made a decision.
Amongst the existing distributed online optimization algorithms, the subgradient descent methods gain considerable attention \cite{MateosD2014Distributed,CaoX2021Decentralized,YanF2013Distributed,CaoX2021DecentrEvent}.
{\color{black}More recently, the authors of \cite{LeeS2016Distributed,LiX2021DistriAgg,ShahrampourS2018Distributed,DixitR2020Online} developed distributed online algorithms on the basis of primal-dual method, gradient tracking, mirror descent approach, and proximal gradient algorithm}, respectively.
The authors of  \cite{YiX2019Distributed} further introduced a primal-dual mirror descent algorithm to address distributed online problems with time-varying coupled inequality constraints over weight-balanced digraphs.
Specifically,  
the algorithms in 
\cite{MateosD2014Distributed,CaoX2021Decentralized,YanF2013Distributed,LeeS2016Distributed,CaoX2021DecentrEvent,LiX2021DistriAgg} 
focused on \textit{static regret}, which can be used to measure the performance of estimating a static target in sensor networks.
While 
the authors of \cite{ShahrampourS2018Distributed,YiX2019Distributed,DixitR2020Online} 
concentrated on \textit{dynamic regret}, {\color{black}which is a more stringent metric than static regret since it allows the best decision in handsignht varies with time.  Therefore, the dynamic regret can reflect the performance of tracking moving targets.}

{\color{black}However,  the aforementioned online algorithms highly rely on the assumption of doubly stochastic weight matrices, which is quite stringent and even deemed impracticable in applications, 
e.g., computer networks intrinsically operate with directionality and it is difficult to construct doubly stochastic weight matrices in a distributed manner. 
Extending distributed algorithms from weight balanced networks to general directed networks is non-trivial \cite{nedic2014distributed,Scutari2018Distributed,CaiK2014Average,Mai2019Distributed,XieP2018Distributed,XiC2017Distributed,PuS2018PushJou,XinR2018A}.} 
%
%
{\color{black}The authors of \cite{Akbari2017Distributed,PangB2019Randomized} proposed distributed online optimization algorithms inspired by the push-sum based algorithm \cite{nedic2014distributed} and the surplus-based method \cite{XiC2017Distributed}, respectively. However, the former is incapable of tackling constrained optimization problems by combining with projected-based methods directly,  
and the latter involves a global parameter depending on weight matrices which should be known a priori. 
}

\begin{table*}[htp]
   \centering
   \begin{threeparttable}[b]
\caption{\textsc{A Comparison of Our Work with The Relevant Distributed Online Optimization Algorithms}}
   \label{Table1}
\begin{tabular}{lcccccc}
\toprule[1pt]
References &  Time-Varying Networks & Unbalancedness &  Constraint Set & Gradient-Free & {\color{black}Dynamic Regret} \\
\midrule
\cite{MateosD2014Distributed}
& \cmark & \xmark &  \xmark  & \xmark & {\color{black}\xmark} \\
\cite{YanF2013Distributed,LeeS2016Distributed}
& \cmark & \xmark &  \cmark  & \xmark & {\color{black}\xmark} \\
\cite{Akbari2017Distributed}
& \cmark & \cmark &  \xmark  & \xmark & {\color{black}\xmark} \\
{\color{black}\cite{CaoX2021Decentralized,CaoX2021DecentrEvent}}
& \xmark & \xmark & \cmark  & \cmark & {\color{black}\xmark} \\
\cite{PangB2019Randomized}
& \xmark & \cmark &  \cmark  & \cmark & {\color{black}\xmark} \\
\cite{Hosseini2016Online}
& \cmark & \cmark & \cmark  & \xmark & {\color{black}\xmark} \\
\cite{ZhangY2019A}
& \xmark & \xmark &  \xmark & \xmark & {\color{black}\cmark} \\
\cite{Yixin2021BanditD}
& \cmark & \xmark &  \cmark  & \cmark & {\color{black}\cmark} \\
\cite{LiX2020Distributed}
& \cmark & \cmark &  \cmark  & \xmark & {\color{black}\cmark} \\
\cite{ShahrampourS2018Distributed}
& \xmark & \xmark & \cmark & \xmark &  {\color{black}\cmark}  \\
{\color{black}\cite{DixitR2020Online,LiX2021DistriAgg}}, \cite{YiX2019Distributed,Luk2019Online}
& \cmark & \xmark  & \cmark & \xmark &   {\color{black}\cmark}    \\
Our work
& \cmark & \cmark & \cmark  & \cmark  &  {\color{black}\cmark}  \\
\toprule[1pt]
   \end{tabular}
  \end{threeparttable}
\end{table*}



{\color{black}In addition, the closed-form expression of cost function or gradient information may be not available in some scenarios, e.g., online source localization, online routing in data networks \cite{Yixin2021BanditD}.  
To relieve this bottleneck, zeroth-order algorithms have gained renewed interests in recent years \cite{Nesterov2017Random,LiuS2018Zeroth,AgarwalA2010Optimal,Shamir2017an,Duchi2015Optimal}, and been studied under distributed setting by combining 
the surplus-based method \cite{PangA2019Randomized},
the primal-dual method \cite{Hajinezhad2019ZONE}, 
the gradient tracking method \cite{Tang2021Distributed}, etc. 
The intuition of these distributed algorithms is  
 constructing gradient estimators from function values and then substituting them for the true gradient under the assumption of doubly stochastic weight matrices or fixed networks. 
Though some efforts have been made to extend zeroth-order distributed algorithms to the online setting \cite{CaoX2021Decentralized,CaoX2021DecentrEvent,PangB2019Randomized}, they can only 
be applicable to fixed networks.}

{\color{black}All the above motivate us to further explore an algorithm that has stronger adaptability.} 
In this work, we study the distributed online optimization problem over time-varying unbalanced
digraphs, under the settings where 
{\color{black}the decisions of nodes are constrained in a convex set and only the local cost function values are revealed to nodes. }
We compare this work with the state-of-the-art on distributed online optimization in Table \ref{Table1}, and
summarize our main contributions as follows:
\begin{enumerate}
\item[(1)]
We propose a novel distributed online constrained optimization algorithm over time-varying digraphs with a local randomized zeroth-order oracle.
    We rigorously analyze the dynamic regret of the algorithm. Our results show that
    the average dynamic regret over a time horizon 
    {\color{black}converges to zero} at a sublinear rate if 
    the accumulated variation grows sublinearly with a specific order.

\item[(2)]  Inspired by the static push-sum protocols developed in \cite{nedic2014distributed,nedic2016Stochastic,Scutari2018Distributed},
    we present a novel methodology that provides a new perspective on the push-sum based distributed online algorithm \cite{Akbari2017Distributed}.
    Specifically,
    by dynamically constructing row-stochastic matrices and rescaling the zeroth-order oracle,
    our algorithm can be applicable to constrained distributed online  optimization problems,
    while the push-sum based online algorithm in \cite{Akbari2017Distributed} is only feasible in the unconstrained case. 

\item[(3)] In comparison with the distributed online optimization algorithms \cite{MateosD2014Distributed,ShahrampourS2018Distributed,YiX2019Distributed,Luk2019Online} relying highly on double stochastic weight matrices, {\color{black}our algorithm can be applicable to unbalanced networks
with column-stochastic weight matrices, and allow nodes to utilize zeroth-order information in lieu of subgradients as well.
    In addition, compared with the fixed networks and static regret considered in \cite{PangB2019Randomized}, 
    our algorithm can be 
    applicable  
    to time-varying networks with dynamic regret guarantees.} 
%
\end{enumerate}

The rest of this paper is organized as follows.
Section \ref{sec-II}
describes the problem of interest.
Section \ref{sec-III} presents the proposed algorithm.
Section \ref{sec-IV} includes our main results on dynamic regret.
Simulations are provided in Section \ref{sec-V}.
We conclude this paper in Section \ref{sec-VI}.

{\bf Notations.}
The superscript `$\text{T}$' denotes vector transposition;
 We use $[T]=\{1,...,T\}$ to denote a set of integers for $T\in\mathbb{N}_{+}$;
The notation $[A]_{ij}$ or $A_{ij}$ denotes the $i,j$-th element of matrix $A$;
$\bm{1}_{n}$ represents a column vector with its all components equaling to one.
We use $\mathcal{P}_{\Omega}(x)$ to denote the projection operation of a point $x$
onto the set $\Omega\subset\mathbb{R}^{m}$, i.e., $\mathcal{P}_{\Omega}(x)=
\arg\min_{w\in\Omega}\Vert w-x\Vert^{2}$.
A subgradient of convex function $f$ at $x$ is represented by $g(x)\in\partial f(x)$, which satisfies
$f(y)\geq f(x)+g^{\text{T}}(x)(y-x)$.  
For two functions $f$ and $h$, the notation $f=\mathcal{O}(h)$ means that there exists a positive constant $\varrho<\infty$ such that $f\leq \varrho h$. 
{\color{black}$\mathbb{B}^{m}=\{x\in\mathbb{R}^{m}|\Vert x\Vert\leq 1\} $ and $\mathbb{S}^{m}=\{x\in\mathbb{R}^{m}|\Vert x\Vert= 1\}$ denote the unit ball and sphere, respectively.}  
We denote by $\mathbb{E}(x)$ the expectation of $x$.

\section{Preliminaries and Problem Formulation}\label{sec-II}



\subsection{Graph Theory}

We consider a network containing $n$ nodes. A sequence of digraphs $\{\mathcal{G}_{t}({\mathcal{V}},\mathcal{E}_{t})\}_{t\in[T]}$ is utilized to
model the interactions among the nodes,
where $\mathcal{V}=\{1,...,n\}$ represents the set of nodes
and $\mathcal{E}_{t}$ denotes the
set of interaction links at $t$. $(i,j)\in\mathcal{E}_{t}$
implies that node $j$ can receive information from node $i$ at $t$.
We denote $\mathcal{N}_{i,t}^{+}=\{j:(j,i)\in\mathcal{E}_{t}\}\cup \{i\}$ and
$\mathcal{N}_{i,t}^{-}=\{j:(i,j)\in\mathcal{E}_{t}\}\cup \{i\}$ as the in-neighbor
and out-neighbor sets of node $i$, respectively.
The digraph $\mathcal{G}_{t}$ is strongly connected if there exists a
path between any pair of distinct nodes.
In addition, $\mathcal{A}_{t}=[a_{ij,t}]\in\mathbb{R}^{n\times n}$ is the weight matrix
induced by $\mathcal{G}_{t}$, 
where $a_{ij,t}>0$ if and only if $(j,i)\in\mathcal{E}_{t}$,
and $a_{ij,t}=0$ otherwise. Moreover, we say matrix $\mathcal{A}_{t}$ is
row-stochastic if $\mathcal{A}_{t}\bm{1}_{n}=\bm{1}_{n}$,
column-stochastic if $\mathcal{A}_{t}^{\text{T}}\bm{1}_{n}=\bm{1}_{n}$,
and doubly stochastic if
$\mathcal{A}_{t}\bm{1}_{n}=\bm{1}_{n}$ and $\mathcal{A}_{t}^{\text{T}}\bm{1}_{n}=\bm{1}_{n}$ hold simultaneously.
The following assumption 
is common in distributed optimization 
\cite{LiX2020Distributed,Akbari2017Distributed,nedic2014distributed}.
\begin{assumption} \label{graph}
The digraph $\mathcal{G}_{t}$, for $ t\in[T]$, satisfies:
\begin{enumerate}
  \item[(a)] \textit{Lower Bounds}: There exists a positive constant $\gamma\in(0,1)$ that lower bounds all nonzero
weights.

  \item[(b)] \textit{Column-Stochasticity}:  $\mathcal{A}_{t}$ is column-stochastic.

  \item[(c)] \textit{Uniformly Jointly Strongly Connected}: The digraphs $\{\mathcal{G}_{t}\}_{t=1}^{\infty}$ are
  uniformly jointly strongly connected, i.e., there exists a positive integer $B$ such that the union
  digraph $\mathcal{G}(\mathcal{V}, \cup_{l=0,...,B-1}\mathcal{E}_{t+l})$ is strongly connected.
\end{enumerate}
\end{assumption}


\subsection{Problem Formulation}


{\color{black}
Consider a network of $n$ nodes, which aim to  
collaboratively solve the following constrained optimization problem over a time horizon 
$T$: 
\begin{eqnarray}\label{Inhand}
\mathop{\text{minimize}}_{x_{t}\in\Omega}~ \sum_{t=0}^{T}F_{t}(x_{t}),\quad F_{t}(x)\triangleq \sum_{i=1}^{n}f_{i,t}(x)
\end{eqnarray}
where $f_{i,t}(x)$ is the local cost function of node $i\in\mathcal{V}$ at round $t\in[T]$. 
We consider the \textit{online} and \textit{bandit feedback} scenario wherein 
each node is only allowed to query 
the values of its local cost function 
after it has made its decision at each round. 
Since the global function $F_{t}$ is not accessible to any node, each node needs to 
interact with its neighbors and make sequential decisions  based on what it ``thinks'' the decisions that the whole  network would make.

The \textit{dynamic regret} 
is commonly adopted in the literature as a performance metric for the online algorithms (see e.g., \cite{Luk2019Online}), 
which is defined as follows for node $j\in\mathcal{V}$ 
\begin{eqnarray}\label{DynamicRegret}
\mathcal{R}_{j}^{d}(T)
=
\sum_{t=0}^{T}F_{t}(x_{j,t})
-\sum_{t=0}^{T}F_{t}(x_{t}^{\star}),
\end{eqnarray}
where $x_{t}^{\star}\in\text{arg}\min_{x\in\Omega}F_{t}(x)$.
The dynamic regret (\ref{DynamicRegret}) measures the difference between the cost incurred by node $j$'s decisions $\{x_{j,t}\}_{t\in[T]}$ against that of a time-varying clairvoyant. 
It is different from that of \cite{LiX2020Distributed,Yixin2021BanditD}, where the decisions of nodes are not required to reach consensus since the global cost function at round $t$ is defined as $F_{t}(X)=\sum_{i=1}^{n}f_{i,t}(x_{i})$ with $X\triangleq col(x_{1},...,x_{n})$ and $X_{t}^{\star}\triangleq col(x_{1,t}^{\star},...,x_{n,t}^{\star})$. 
In contrast to the \textit{static regret} \cite{MateosD2014Distributed} with fixed benchmark $x^{\star}=\text{arg}\min_{x\in\Omega}\sum_{t=0}^{T}F_{t}(x)$, the dynamic regret (\ref{DynamicRegret}) allows the best decision varies with time, and thus can encompass the former as a special case. 
However, node $i\in\mathcal{V}$ can not access to $f_{i,t}$ when making decision $x_{i,t}$, 
which implies that the online algorithm is not able to track $x_{t}^{\star}$ well if $x_{t}^{\star}$ deviates from the past significantly. 
Therefore, 
a regularity measure named ``path variation'' \cite{Luk2019Online} needs to be introduced 
to reflect the changes of successive minimizers  
%
%
%
%
%
%
\begin{eqnarray}\label{Curvari}
\mathcal{C}_{T}
=
\sum_{t=0}^{T}\Vert x_{t+1}^{\star}-x_{t}^{\star}\Vert.
\end{eqnarray}

{\color{black}
The main objective of this work is to design a 
zeroth-order distributed online optimization algorithm over $\mathcal{G}_{t}$ 
such that 
the dynamic regret (\ref{DynamicRegret}) is upper bounded sublinearly with respect to $T$.
We impose the following standard assumption in (distributed) bandit online optimization \cite{AgarwalA2010Optimal,Yixin2021BanditD,YanF2013Distributed}.} 
\begin{assumption}\label{ConvexandBoundedG}
For $ i\in\mathcal{V}$ and $t\in[T]$,  it satisfies
\begin{enumerate}
  \item[(a)] $\Omega\subseteq\mathbb{R}^{m}$ is a non-empty, convex and closed set.
  {\color{black}Moreover, there exist $r, R>0$ such that 
  \begin{eqnarray}
{\color{black} r\mathbb{B}^{m}\subseteq \Omega\subseteq R\mathbb{B}^{m},}
\end{eqnarray}
and $r$ is known a priopri.}
\item[(b)] The subgradient of $f_{i,t}$ is bounded over $\Omega$, i.e., there exists a constant $G>0$ such that $\Vert g_{i,t}(x)\Vert\leq G$ holds for $ x\in\Omega$,  where $g_{i,t}(x)\in \partial f_{i,t}(x)$.
\end{enumerate}
\end{assumption}
}

\section{Algorithm Development}\label{sec-III}


{\color{black}
To develop the zeroth-order algorithm, we first define a smoothed version of $f_{i,t}(x)$, $i\in\mathcal{V}$, as follows
\begin{eqnarray}
{\color{black}\hat{f}_{i,t}(x)\triangleq\mathbb{E}_{\zeta\in\mathbb{B}^{m}}[f_{i,t}(x+\mu\zeta_{i,t})], ~~x\in (1-\xi)\Omega, }
\end{eqnarray}
where 
$\zeta_{i,t}$ is a vector selected uniformly at random from $\mathbb{B}^{m}$, 
 $\xi\in(0,1)$ is a shrinkage parameter, 
$\mu\in(0, r\xi]$,
and $(1-\xi)\Omega$ is shorthand for $\{(1-\xi)x: x\in\Omega\}$.  
Under these settings, it can be seen that $x+\mu\zeta_{i,t}\in\Omega$ for 
 any $x\in(1-\xi)\Omega$ and $\zeta_{i,t}\in\mathbb{S}^{m}$.
%
In this paper, we adopt the following local zeroth-order oracle to avoid explicit subgradient calculations 
\begin{eqnarray}\label{GForacle}
\hat{g}_{i,t}(x)=\frac{m}{\mu}\left(f_{i,t}(x+\mu\zeta_{i,t})-f_{i,t}(x)\right)\zeta_{i,t}, 
\end{eqnarray}
where $x\in(1-\xi)\Omega$, $\zeta_{i,t}\in\mathbb{S}^{m}$.
In fact, $\hat{g}_{i,t}(x)$ is an unbiased gradient estimator of $\hat{f}_{i,t}(x)$. Notably, Gaussian random variables are used to construct zeroth-order oracles in \cite{PangB2019Randomized,PangA2019Randomized,Yuan2014Randomized}, which cannot be applied in our setting. The reason is that Gaussian random variables do not have finite support such that the perturbation $x+\mu\zeta_{i,t}$ may lie outside of $\Omega$.} 

In our algorithm, each node $i\in\mathcal{V}$ maintains a vector $x_{i,t}\in\mathbb{R}^{m}$ as well as a scalar  $\phi_{i,t}\in\mathbb{R}$ at per round $t\in\mathbb{N}$, which are initialized with $x_{i,0}\in(1-\xi)\Omega$ and $\phi_{i,0}=1$, respectively. 
At per round $t$, each node $i$ sends its decision $x_{i,t}$ and the scalar $\phi_{i,t}$ to its out-neighbors on the
basis of $\mathcal{G}_{t}$, and then performs the following updates: 
%
%
%
%
\begin{eqnarray}\label{Algorithm_1}
\phi_{i,t+1}
&\hspace*{-0.5em}=\hspace*{-0.5em}&
\sum_{j=1}^{n}a_{ij,t}\phi_{j,t},  \label{Alg_1} \\
x_{i,t+1}
&\hspace*{-0.5em}=\hspace*{-0.5em}&
\mathcal{P}_{{\color{black}(1-\xi)\Omega}}\left(\sum_{j=1}^{n}b_{ij,t}x_{j,t}-
\frac{1}{\phi_{i,t+1}}\alpha_{t} {\color{black}\hat{g}_{i,t}(x_{i,t})}\right),   \label{Alg_2}
\end{eqnarray}
where
$\alpha_{t}$ is the step size that will be specified later, 
$b_{ij,t}{\color{black}\triangleq}\frac{1}{\phi_{i,t+1}}a_{ij,t}\phi_{j,t}$ is the entry at $i$-th row and $j$-th column of matrix $\mathcal{B}_{t}$. 
Under Assumption \ref{graph}(b),
it is straightforward to verify that $\mathcal{B}_{t}$, $\forall t\geq 0$, now becomes a  row stochastic matrix
{\color{black}by recalling (\ref{Alg_1}) and the definition of $b_{ij,t}$. 
By this means, the term $\sum_{j=1}^{n}b_{ij,t}x_{j,t}$ in (\ref{Alg_2}) is designed for consensus since our setting requires that each node can make sequential decisions  to minimize the dynamic regret.
The intuition of $\phi_{i,t+1}$ in (\ref{Alg_1}) is that it is designed to rescale the zeroth-order oracle $\hat{g}_{i,t}(x_{i,t})$ in (\ref{Alg_2}), which aims to address the unbalancedness issue of time-varying digraphs. 
In addition, the projection operation  $\mathcal{P}_{(1-\xi)\Omega}(\cdot)$ is used to guarantee 
the feasibility of the quarry points, i.e.,  
 $x_{i,t}+\mu\zeta_{i,t}\in\Omega$ always holds for all $i\in\mathcal{V}$ and $t\in[T]$.} 
Our algorithm is fully distributed that each node only requires the knowledge of itself as well as information from its immediate neighbors to carry out the above updates.
The proposed algorithm is summarized in Algorithm \ref{algOnline}.

\begin{algorithm}[t]
\caption{{\color{black}Zeroth-order Distributed Online Optimization Algorithm} ~~{\color{black}---from the view of node $i$}}
\label{algOnline}
\begin{algorithmic}[1]
    \State \textbf{Initialize:} Set $x_{i,0}\in{\color{black}(1-\xi)\Omega}$ and $\phi_{i,0}=1$, $ i\in\mathcal{V}$. 

    \For {$t=0,1,2,...,T$}

    \State Broadcast $x_{i,t}$ and $\phi_{i,t}$ to nodes $j\in\mathcal{N}_{i,t}^{-}$.

    \State Calculate {\color{black}$\hat{g}_{i,t}(x_{i,t})$} via (\ref{GForacle}).
    \State Update  auxiliary variable $\phi_{i,t}$ via (\ref{Alg_1}).
    \State Update node's decision $x_{i,t}$ via (\ref{Alg_2}).
    \EndFor


\end{algorithmic}
\end{algorithm}


Several features distinguish our algorithm from the relevant distributed (online) optimization methods over unbalanced digraphs. We summarize their connections as follows:

\begin{enumerate} 
 {\color{black}
  \item[(a)] The $\mathcal{AB}$/push-pull algorithms \cite{XinR2018A,PuS2018PushJou} and our algorithm all utilize row stochastic and column stochastic weight matrices simultaneously. Different from \cite{XinR2018A,PuS2018PushJou} that the weight matrices are constructed independently, our algorithm only invokes a sequence of column stochastic matrices $\mathcal{A}_{t}$, and we obtain $\mathcal{B}_{t}$ by performing the designated transformation on $\mathcal{A}_{t}$.
  
   \item[(b)]  The push-sum based distributed online optimization algorithm in \cite{Akbari2017Distributed} will be invalid when it  comes to problems with constraint set. 
  The reason is that the balance properties guaranteed by taking specific ratios will be violated if we incorporate the projection operation into push-sum algorithms \cite{nedic2014distributed,nedic2016Stochastic,Scutari2018Distributed} directly. 
   In contrast, our algorithm can deal with set constraint and be applicable to scenarios where the explicit subgradients of the cost functions are not available. 
   
    \item[(c)] 
    The static distributed algorithm in \cite{Mai2019Distributed} also adjusts the gradient by a designed vector, 
%
    which is essentially the estimate of left 
      eigenvector of the weight  matrix. However, it is not clear how to extend the algorithm to time-varying digraphs, since it is unrealistic to estimate the time-varying Perron vectors in a distributed manner.    
  Our algorithm can be applicable to
      time-varying digraphs. 
      Furthermore, each node $i$ transmits the scalar $\phi_{i,t}$ at $t$ instead of a vector, which considerably lightens the  communication burden when compared with \cite{Mai2019Distributed}.
   
 }
\end{enumerate} 

\section{Convergence Analysis}\label{sec-IV}

In this section, we provide the theoretical analysis for the proposed algorithm.
The main challenge in the convergence analysis lies in the combined effects
of the projection operators, the time-varying unbalanced networks, and the embedded
randomized zeroth-order oracles.
 In addition, the considered dynamic regret makes our analysis more challenging than that of the static regret scenario, since the cost
functions and the minimizers are allowed to drift over time simultaneously in our setting. 
We begin with a few preliminary results,
followed by the detailed analysis on dynamic regret of the proposed algorithm.


\subsection{Preliminary results}

%

{\color{black}Let $\mathcal{F}_{t}$ denote the $\sigma$-field generated by the entire history of the random variables 
up to round $t$, i.e., $\mathcal{F}_{t}\triangleq \sigma(\zeta_{1,0},...,\zeta_{n,0},...,\zeta_{1,t},...,\zeta_{n,t})$.}   
The following Lemma collected from \cite{Yixin2021BanditD} outlines  some important properties of 
{\color{black}$\hat{f}_{i,t}(x)$ and $\hat{g}_{i,t}(x)$.} 
{\color{black}
\begin{lemma}\cite{Yixin2021BanditD}\label{gradient-freePro}
Suppose Assumption 2 holds. For $ i\in\mathcal{V}$, and $t\geq 0$, the following properties hold: 
\begin{enumerate}
  \item[(a)] $\hat{f}_{i,t}(x)$ is convex and $G$-Lipschitz continuous on $(1-\xi)\Omega$. Moreover, 
    \begin{eqnarray}
  f_{i,t}(x)
\leq
  \hat{f}_{i,t}(x)
\leq
  f_{i,t}(x)+\mu_{i} G.
  \end{eqnarray} 
  for any $x\in(1-\xi)\Omega$. 
%

  \item[(b)] $\hat{f}_{i,t}(x)$ is differentiable on $(1-\xi)\Omega$ even if ${f}_{i,t}(x)$ is not, and it satisfies
  \begin{eqnarray}
    \nabla \hat{f}_{i,t}(x)=\mathbb{E}_{\zeta_{i,t}\in\mathbb{S}^{m}}[\hat{g}_{i,t}(x)], ~x\in(1-\xi)\Omega. 
  \end{eqnarray}
  Moreover, $\nabla \hat{f}_{i,t}(x)$ is Lipschitz continuous on $(1-\xi)\Omega$ with $L_{i}\triangleq mG/\mu_{i}$. 
%

  \item[(c)] $\Vert\hat{g}_{i,t}(x)\Vert\leq mG$ holds for any $x\in(1-\xi)\Omega$. 
%
\end{enumerate}
\end{lemma}
}
{\color{black}Recalling that $\mathcal{B}_{t}$ is a row stochastic matrix for any $ t\geq 0$, and it thus enjoys the following properties.}
\begin{lemma}\cite{XieP2018Distributed}\label{Convergence-z}
Suppose Assumptions \ref{graph} holds.
For $s\geq t$, define $\mathcal{B}(s:t):=\mathcal{B}_{s-1}\cdots\mathcal{B}_{t}$ with the convention
$\mathcal{B}(t:t)=I$ and $b_{ij}(s:t)$ being the entries of $\mathcal{B}(s:t)$.
Then, there is a sequence of normalized vectors $\{\pi_{t}\}_{t\geq 0}$ with
$\pi_{t}=[\pi_{1,t},...,\pi_{n,t}]^{\text{T}}$ and $\bm{1}_{n}^{\text{T}}\pi_{t}=1$, such that
\begin{enumerate}
  \item[(a)] There exists constants $C>0$ and $\lambda\in(0,1)$, such that $|b_{ij}(s:t)-\pi_{j,t}|\leq C\lambda^{s-t}$ holds
  for 
  $ i,j\in\mathcal{V}$;

  \item[(b)] There exists a constant $\beta\geq b_{\min}^{(n-1)B}$ such that $\pi_{i,t}\geq \beta $ for $ i\in\mathcal{V}$ and $t\geq 0$, where $b_{\min}=\min_{i,j\in\mathcal{V},t\geq 0}\{b_{ij,t}|j\in\mathcal{N}_{i,t}^{+}\cup \{i\}\}$;

  \item[(c)] $\pi_{t}^{\text{T}}=\pi_{t+1}^{\text{T}}\mathcal{B}_{t}$.
\end{enumerate}
\end{lemma}

The following lemma is collected from 
\cite{Scutari2018Distributed}, which reveals that the scalar $\phi_{i,t}$
is bounded for $ i\in\mathcal{V}$, $t\geq 0$.
\begin{lemma}\cite{Scutari2018Distributed}\label{BoundLe}
Suppose Assumption \ref{graph} holds.
Let $\{\phi_{i,t}\}_{t\geq 0}$, $\forall i\in\mathcal{V}$, be the sequence generated by (\ref{Alg_1}). Then, there exists  constants $\theta >0$ and $\varpi >0$ such that
\begin{eqnarray}\label{BoundPhi}
 \theta^{-1} \leq \phi_{i,t} \leq \varpi.
\end{eqnarray}
\end{lemma}

\begin{remark}
In fact, $\theta$ and $\varpi$ in Lemma \ref{BoundLe} can be precisely represented by $\gamma^{-2(n-1)B}$ and $n-\gamma^{2(n-1)B}$, respectively, where $\gamma$ and $B$ are constants defined in Assumption \ref{graph}.
\end{remark}

Furthermore,
the following lemma establishes a relationship between $\phi_{i,t+1}$ and $\pi_{i,t+1}$, which plays a critical role
in the dynamic regret analysis in the sequel.
\begin{lemma}\label{Rescale}
Let $\{\phi_{i,t+1}\}_{t\geq 0}$, $\forall i\in\mathcal{V}$, be the sequence generated by (\ref{Alg_1}), and $\pi_{i,t+1}$
be the $i$-th component of $\pi_{t+1}$ clarified in Lemma \ref{Convergence-z}. Then, for a given time horizon $T$, the following inequality holds
\begin{eqnarray}\label{RescalePhi}
\left|\frac{1}{n}\phi_{i,t+1}-\pi_{i,t+1}\right|
&\hspace*{-0.5em}\leq\hspace*{-0.5em}&
C\lambda^{T-t-1}.
\end{eqnarray}
\end{lemma}

\textit{Proof:}
See Appendix A.
$\hfill \blacksquare$
\subsection{{\color{black}Dynamic} Regret Analysis}

As detailed in the previous sections,
the combination of projection operation and time-varying unbalanced networks brings new
challenges in the dynamic regret analysis compared to the existing works. 
Particularly, we cannot
concentrate on the evolution of the average process 
$\frac{1}{n}\sum_{i=1}^{n}x_{i,t}$ 
as many distributed algorithms \cite{Akbari2017Distributed,nedic2014distributed,Luk2019Online} did.
{\color{black}To cope with this challenge, we  construct time-varying row-stochastic matrices $\mathcal{B}_{t}$, which enables us to turn attention from the average process to an auxiliary vector defined as
\begin{eqnarray}\label{x_bar-1} 
\bar{x}_{t}\triangleq\sum_{i=1}^{n}\pi_{i,t}x_{i,t}
\end{eqnarray}
with $\pi_{t}=[\pi_{1,t},...,\pi_{n,t}]$ defined in Lemma \ref{Convergence-z}.}
The following lemma establishes an upper bound on the expected disagreement between 
$\bar{x}_{t}$ and $x_{i,t}$ for $i\in\mathcal{V}$. 

\begin{lemma}\label{theoremcons}
Suppose Assumptions 1-2 hold.
Let $\{x_{i,t}\}_{t\in[T]}$ be the sequence generated by Algorithm \ref{algOnline}. 
Then, for $i\in\mathcal{V}$, $t\in[T]$, 
%
%
\begin{eqnarray} \label{Theorem-1}
&&\mathbb{E}\left[\Vert x_{i,t+1}-\bar{x}_{t+1}\Vert\right]  \notag\\
&&\leq
n{\color{black}m}\theta CG\sum_{l=1}^{t+1}\lambda^{t+1-l}\alpha_{l-1} +C\lambda^{t+1}\sum_{j=1}^{n}\Vert x_{j,0}\Vert. 
\end{eqnarray}
\end{lemma}

\textit{Proof:}
See Appendix B. 
$\hfill \blacksquare$

\begin{remark}
{\color{black}
Note that the first term on the right side of (\ref{Theorem-1}) can be calculated as follows 
\begin{eqnarray} \label{remar1}
&&\hspace{-1em}\sum_{l=1}^{t+1}\lambda^{t+1-l}\alpha_{l-1}
=
\sum_{l=1}^{\lfloor t/2\rfloor}\lambda^{t-l+1}\alpha_{l-1}
+\sum_{l=\lfloor t/2\rfloor +1}^{t+1}\lambda^{t-l+1}\alpha_{l-1}  \notag\\
&&\hspace{-1em}=
\lambda^{\lceil t/2\rceil}\sum_{l=1}^{\lfloor t/2\rfloor}\lambda^{\lfloor t/2\rfloor-l+1}\alpha_{l-1}
+\sum_{l=\lfloor t/2\rfloor +1}^{t+1}\lambda^{t-l+1}\alpha_{l-1} \notag\\
&&\hspace{-1em}\leq
\xi^{\lceil t/2\rceil}\frac{\lambda}{1-\lambda}\sup_{t\geq 0}\alpha_{l}
+\frac{\lambda}{1-\lambda}\sup_{l\geq\lfloor t/2\rfloor}\alpha_{l},   \notag
\end{eqnarray}  
where $\lfloor \cdot\rfloor$ and $\lceil \cdot\rceil$ denote the floor function and ceiling function, respectively. 
If the step size $\alpha_{t}$ tends to zero, e.g., $\alpha_{t}=\mathcal{O}(\frac{1}{\sqrt{t+1}})$, 
then the two terms in the above inequality both tend to zero since $\xi^{\lceil t/2\rceil}\to 0$ and $\sup_{l\geq\lfloor t/2\rfloor}\alpha_{l}\to 0$ when $t\to\infty$.  
Thus we can obtain $\lim_{t\to\infty}\mathbb{E}[\lVert x_{i,t+1}-x_{j,t+1} \rVert ] =0$, 
which implies that all nodes can achieve consensus in expectation when $t\to\infty$. 
}
\end{remark}
To obtain the desired regret bound, it suffices to
derive an upper bound on the expected optimality gap. {\color{black}However, $x_{t}^{\star}$ may lie outside of $(1-\xi)\Omega$, which bring new challenge to our analysis. To resolve this issue, we first establish a relation on $\mathbb{E}[f_{i,t}(\bar{x}_{t})-f_{i,t}((1-\xi)x_{t}^{\star})]$,  and then the desired result can be obtained by bounding the difference between $f_{i,t}((1-\xi)x_{t}^{\star})-f_{i,t}(x_{t}^{\star})$.} The following result provides a bound on the weighted sum of the terms $\mathbb{E}[f_{i,t}(\bar{x}_{t})-f_{i,t}(x_{t}^{\star})]$, which is pivotal to our dynamic regret analysis.

\begin{theorem}\label{Lemma_4}
Suppose Assumptions 1-2 hold. Let $\{x_{i,t}\}_{t\in[T]}$ be the sequence generated by Algorithm 1 with positive and non-increasing step size $\alpha_{t}$.
Then,
{\color{black}\begin{eqnarray}\label{Theorem22}
&&\hspace{-1.5em}\sum_{t=0}^{T}\sum_{i=1}^{n}\frac{\pi_{i,t+1}}{\phi_{i,t+1}}\mathbb{E}\left[f_{i,t}(\bar{x}_{t})-f_{i,t}(x_{t}^{\star})\right]   \notag\\
&&\hspace{-1.5em}\leq
2nm\theta^{2}RLCG\sum_{t=0}^{T}\sum_{l=1}^{t}\lambda^{t-l}\alpha_{l-1}
+\frac{m^{2}\theta^{2}G^{2}}{2}\sum_{t=0}^{T}\alpha_{t} \notag\\
&&\hspace{-1.5em}~~~+2nm^{2}\theta^{2}CG^{2}\sum_{t=0}^{T}\sum_{l=1}^{t+1}\lambda^{t+1-l}\alpha_{l-1}  
 +\frac{2R^{2}}{\alpha_{T}} \notag\\
&&\hspace{-1.5em}~~~+2\left(RL+mG\lambda\right)\theta C\sum_{j=1}^{n}\Vert x_{j,0}\Vert\sum_{t=0}^{T}\lambda^{t} 
+\frac{2R}{\alpha_{T}}\mathcal{C}_{T}  \notag\\
&&\hspace{-1.5em}~~~+\mu\theta G (T+1) +\theta \xi GR(T+1). 
\end{eqnarray}}
%
%
%
%
%
%
\end{theorem}

\textit{Proof:}
See Appendix C.
$\hfill \blacksquare$

Albeit Theorem \ref{Lemma_4} has established a relation on the term $\mathbb{E}[f_{i,t}(\bar{x}_{t})-f_{i,t}(x_{t}^{\star})]$, it should be noted here that
this term is scaled by $\frac{\pi_{i,t+1}}{\phi_{i,t+1}}$ and we cannot determine the sign of the difference (i.e., $f_{i,t}(\bar{x}_{t})-f_{i,t}(x_{t}^{\star})$) for particular $i\in\mathcal{V}$ since $x_{t}^{\star}$ is the minimizer of the global cost function rather than each local cost function, which bring new challenge in deriving the dynamic regret.
{\color{black}Fortunately, 
we will soon find that the scalar $\frac{1}{n}\phi_{i,t+1}$ plays an important role in counteracting the impact of unbalancedness,    
and hence we can establish an upper bound of the dynamic regret (\ref{DynamicRegret}) in the following theorem.

}

{\color{black}
\begin{theorem}[Dynamic Regret Bound]\label{NetworkRegret}
Suppose Assumptions 1-2 hold. For a given time horizon $T$,  set $\mu=\frac{r}{\sqrt{T+1}}$ and $\xi=\frac{1}{\sqrt{T+1}}$. 
Let $\{x_{i,t}\}_{t\in[T]}$ be the sequence generated by Algorithm \ref{algOnline} with  $\alpha_{t}=\frac{1}{m\sqrt{t+1}}$. 
Then, 
\begin{eqnarray}\label{NetRegRes}
\hspace{-1em}\mathcal{R}_{j}^{d}(T)
&\hspace*{-0.5em}\leq\hspace*{-0.5em}&
\mathcal{K}_{1}+\mathcal{K}_{2}\sqrt{T+1}
\end{eqnarray}
where
\begin{eqnarray}
\mathcal{K}_{1}
&\hspace*{-0.5em}=\hspace*{-0.5em}&
\frac{nC(2RL\theta +2m\theta\lambda G+G)}{1-\lambda}\sum_{i=1}^{n}\Vert x_{i,0}\Vert 
+\frac{2n^{2}\theta CGR}{\lambda(1-\lambda)},  \notag\\
\mathcal{K}_{2}
&\hspace*{-0.5em}=\hspace*{-0.5em}&
\frac{2n^{2}\theta CG}{1-\lambda}\left(2\theta RL+2m\theta G+G \right)
+nm\theta^{2}G^{2}  \notag\\
&&+2nmR^{2}+nr\theta G+n\theta GR+2nmR\mathcal{C}_{T}. 
%
\end{eqnarray}
\end{theorem}}

\textit{Proof:}
See Appendix D.
$\hfill \blacksquare$

{\color{black}
Theorem \ref{NetworkRegret} shows that if $\mathcal{C}_{T}$ is bounded, then the average dynamic regret over $T$  asymptotically converges zero at the rate of $\mathcal{O}(T^{-\frac{1}{2}})$, which 
matches the best static regret that can be achieved in literature for convex cost functions \cite{MateosD2014Distributed,CaoX2021Decentralized,CaoX2021DecentrEvent,YanF2013Distributed,Akbari2017Distributed,Hosseini2016Online}.
In fact, 
$\mathcal{C}_{T}$ is not necessary to be bounded. 
If $\mathcal{C}_{T}$ increases sublinearly with the rate ranged from
zero to $\mathcal{O}(T^{\frac{1}{2}})$, then we have $\lim_{T\to\infty}\frac{\mathcal{C}_{T}\sqrt{T+1}}{T}=0$, which implies that
the average dynamic regret over $T$ asymptotically converges to zero at a sublinear rate. 
}
{\color{black}\begin{remark}
Zeroth order optimization algorithms usually suffer from deterioration in convergence as the problem dimension increases, which is a typical limitation of these algorithms in high dimension problems. 
From Lemma \ref{gradient-freePro}(c), it is readily seen that the upper bound of $\Vert \hat{g}_{i,t}(x)\Vert$ depends on the dimension $m$, 
which is essentially the
penalty incurred by the use of zeroth-order oracle instead of the real subgradient. 
The dimension dependency of our algorithm is $\mathcal{O}(m)$, which is identical to that of \cite{TangY2020zeroth}, and better than 
$\mathcal{O}(m^{2})$ in \cite{CaoX2021Decentralized,CaoX2021DecentrEvent,PangB2019Randomized,Yuan2014Randomized,Nesterov2017Random,AgarwalA2010Optimal}. 
The optimal dimension dependency $\mathcal{O}(\sqrt{m})$ is obtained in \cite{Shamir2017an}. However, the algorithm in \cite{Shamir2017an} is centralized, and the author only considered static regret rather than dynamic regret. 
\end{remark}}

\subsection{Discussion}

In this section, we discuss about the counterpart of Algorithm 1 when subgradients of local cost functions are available.
In this scenario,
a new algorithm can be directly obtained by utilizing the subgradient in lieu of the
local randomized zeroth-order oracle {\color{black}$\hat{g}_{i,t}(x_{i,t})$} in (\ref{Alg_2}). Then, each node $i\in\mathcal{V}$ 
performs the following updates\footnote{In this paper, the notation of $x_{i,t}$ and $\phi_{i,t}$, will be kept the same across different algorithms to avoid notational clutter, and it is clear from context which method is in question.}:
\begin{eqnarray}\label{Algorithm_2}
\phi_{i,t+1}
&\hspace*{-0.5em}=\hspace*{-0.5em}&
\sum_{j=1}^{n}a_{ij,t}\phi_{j,t},  \label{Alg_Gra11} \\
x_{i,t+1}
&\hspace*{-0.5em}=\hspace*{-0.5em}&
\mathcal{P}_{\Omega}\left(\sum_{j=1}^{n}b_{ij,t}x_{j,t}-\frac{1}{\phi_{i,t+1}}\alpha_{t}g_{i,t}(x_{i,t})\right),  \label{Alg_Gra12}
\end{eqnarray}
where $g_{i,t}(x_{i,t})$ denotes the subgradient of $f_{i,t}(x)$  evaluated at  $x_{i,t}$. This algorithm is summarized in Algorithm \ref{algOnlineGradient}.
\begin{algorithm}[t]
\caption{Distributed Online Optimization Algorithm ~{\color{black}---from the view of node $i$}}
\label{algOnlineGradient}
\begin{algorithmic}[1]

    \State \textbf{Input:} time horizon $T$; for all $ i\in\mathcal{V}$, set $x_{i,0}\in\Omega$ and $\phi_{i,0}=1$.

    \For {$t=0,1,2,...,T$}

    \State Broadcast $x_{i,t}$ and $\phi_{i,t}$ to nodes $j\in\mathcal{N}_{i,t}^{-}$.

    \State Update  auxiliary variable $\phi_{i,t}$ via (\ref{Alg_Gra11}).
    \State Update node's decision $x_{i,t}$ via (\ref{Alg_Gra12}).

    \EndFor
    \State \textbf{Output}: $\{x_{i,t}\}_{i=1}^{n}$.

\end{algorithmic}
\end{algorithm}

Accordingly, the determined form of the dynamic regret (\ref{DynamicRegret}) can be represented by $\tilde{\mathcal{R}}_{j}^{d}(T)$ as follows
\begin{eqnarray}\label{AnoDynamicRegret}
\tilde{\mathcal{R}}_{j}^{d}(T)
=
\sum_{t=0}^{T}f_{t}(x_{j,t})
-\sum_{t=0}^{T}f_{t}(x_{t}^{\star}),
\end{eqnarray}
Now we provide an upper bound of the dynamic regret (\ref{AnoDynamicRegret}) under Algorithm 2. The following theorem reveals that the average dynamic regret over time horizon $T$ asymptotically converges to
zero at a sublinear rate for convex cost functions provided that the accumulated variation grows sublinearly with the rate ranged from zero to $\mathcal{O}(T^{\frac{1}{2}})$,  which reflects that the convergence of Algorithm \ref{algOnline} is essentially not affected by the incorporated zeroth-order oracle.

\begin{theorem}\label{Theo4G}
Suppose Assumptions 1-2 hold. Let $\{x_{i,t}\}_{t\geq 0}$ be the sequence generated by Algorithm \ref{algOnlineGradient} with  $\alpha_{t}=\frac{1}{\sqrt{t+1}}$.
Then, for a given time horizon $T$, the dynamic regret defined in (\ref{AnoDynamicRegret}) satisfies
\begin{eqnarray}\label{NetReg}
\tilde{\mathcal{R}}_{j}^{d}(T)
&\hspace*{-0.5em}\leq\hspace*{-0.5em}&
\mathcal{K}_{3}+\mathcal{K}_{4}\sqrt{T},
\end{eqnarray}
where
\begin{eqnarray}
\mathcal{K}_{3}
&\hspace*{-0.5em}=\hspace*{-0.5em}&
\frac{(2\theta \lambda+\theta +2)nCG}{1-\lambda}\sum_{i=1}^{n}\Vert x_{i,0}\Vert
+\frac{2n^{2}\theta CGR}{\lambda(1-\lambda)},  \notag\\
\mathcal{K}_{4}
&\hspace*{-0.5em}=\hspace*{-0.5em}&
\frac{2(3\theta+2)n^{2}\theta CG^{2}}{1-\lambda}
+n\theta^{2}G^{2}
+2nR\mathcal{C}_{T}
+2nR^{2}.  \notag
\end{eqnarray}
\end{theorem}

\textit{Proof:}
The proof follows the similar line of Theorem 1-2, we thus omit it.
$\hfill \blacksquare$

%
%
%
%
%


\section{Numerical Examples}\label{sec-V}

In this section, we numerically demonstrate the performance of the proposed algorithms. First, motivated by \cite{Luk2019Online}, we validate our theoretical findings by a numerical example. After that, we apply the proposed algorithms to a distributed target tracking problem, which has been widely investigated in literature, e.g., \cite{ShahrampourS2018Distributed,ZhangY2019A}. 
 Finally, we investigate the dynamic sparse signal recovery problem \cite{DixitR2020Online} and compare our algorithm with the existing ones. 

\subsection{Example I}
\begin{figure}[ptbh]
\centerline{
\includegraphics[width=6cm]{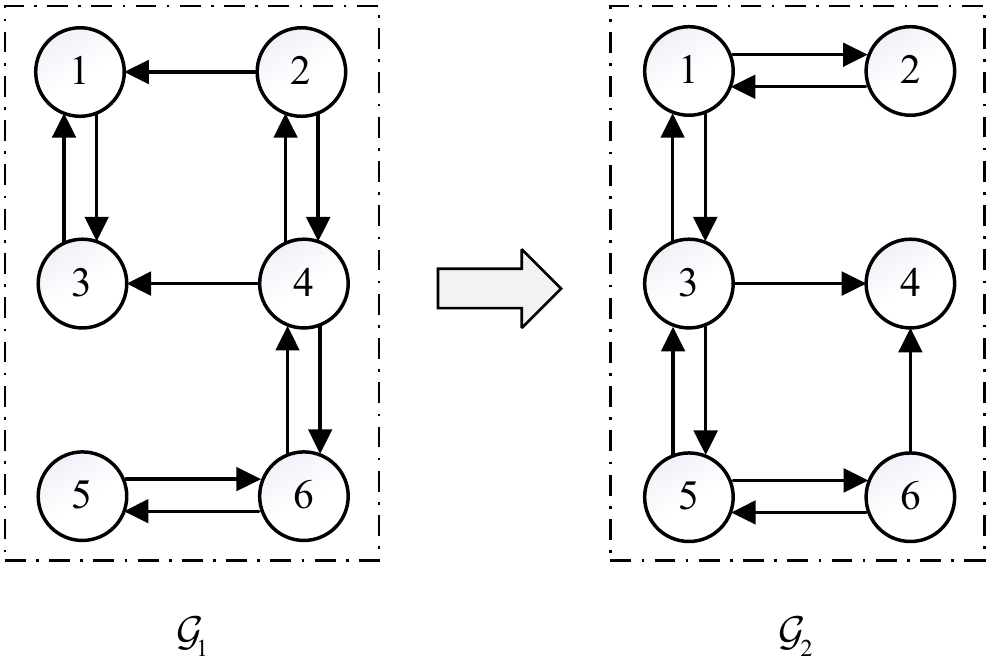} \vspace{-0ex}}
\caption{The time-varying unbalanced digraphs.}
\label{figurea}
\end{figure}

\begin{figure}[ptbh]
\centerline{
\includegraphics[width=9cm]{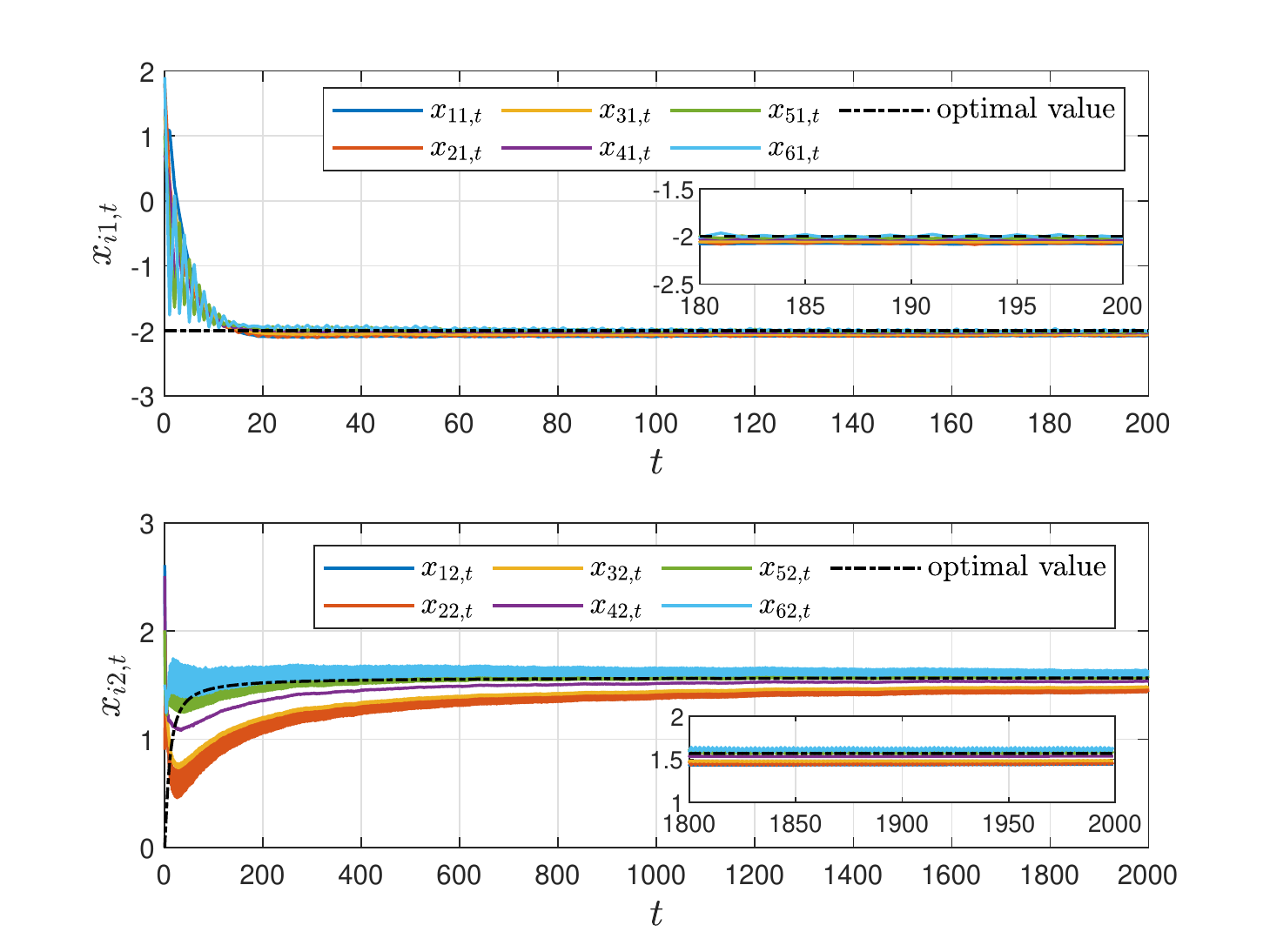} \vspace{-0ex}}
\caption{Trajectories of the optimal solution $x_{t}^{\star}$ and the decisions $x_{i,t}$, $i\in\mathcal{V}$, generated by Algorithm 1. }
\label{figureexa11}
\end{figure}

Consider a group of six agents modeled by the time-varying digraphs $\mathcal{G}_{t}$ depicted in Fig. \ref{figurea}.
The dynamic local cost function at time $t$ is given by $f_{i,t}(\bm{x})=\frac{i}{84}x_{1}^{4}+\frac{i-1}{15}(x_{1}^{2}+x_{2}^{2})+\frac{2i+1}{4}x_{1}-\frac{2(i-3)}{3}\hbar(t)x_{2}$, $i\in\mathcal{V}$ and $\bm{x}=(x_{1},x_{2})^{\text{T}}$ is constrained in a box set $\Omega=\{-3\leq x_{1}\leq 2, 0\leq x_{2}\leq 3\}$.
Thus, the global cost function at time $t$ can be described as $f_{t}(\bm{x})=\frac{1}{4}x_{1}^{4}+x_{1}^{2}+x_{2}^{2}+12x_{1}-2\hbar (t)x_{2}$ with $\bm{x}\in\Omega$.
We let $\hbar (t)=\arctan(t/10)$. It is easy to verify that $\bm{x}_{t}^{\star}=(-2,\arctan(t/10))$, $t\geq 0$. In our simulations, we set $\mu=10^{-3}$ and $\xi=0.02$. 
The graph $\mathcal{G}_{t}$ changes by the order $\mathcal{G}_{1}\to\mathcal{G}_{2}\to\mathcal{G}_{1}\to\cdots$ throughout the whole process.
For the weight matrix $\mathcal{A}_{t}$, we set
$a_{ij,t}=1/|\mathcal{N}_{j,t}^{-}|$ if $a_{ij,t}>0$, 
and $a_{ii,t}=1-\sum_{j=1}^{n}a_{ji,t}$, 
where $|\mathcal{N}_{j,t}^{-}|$ denotes the number of elements in $\mathcal{N}_{j,t}^{-}$.
Under these settings, all assumptions in this paper evidently hold.
We run Algorithm 1 with  $\alpha_{t}=\frac{1}{2\sqrt{t+1}}$.
The trajectories of the optimal solution $x_{t}^{\star}$ and the decision variables $x_{i,t}$, $i\in\mathcal{V}$,  are shown in Fig. \ref{figureexa11} with $x_{ij,t}$ being the $j$-th entry of the decision variable made by agent $i$ at round $t$.
It can be seen that all agents's decision variables approach the optimal solution $x_{t}^{\star}$.

\begin{figure}[ptbh]
\centerline{
\includegraphics[width=8cm]{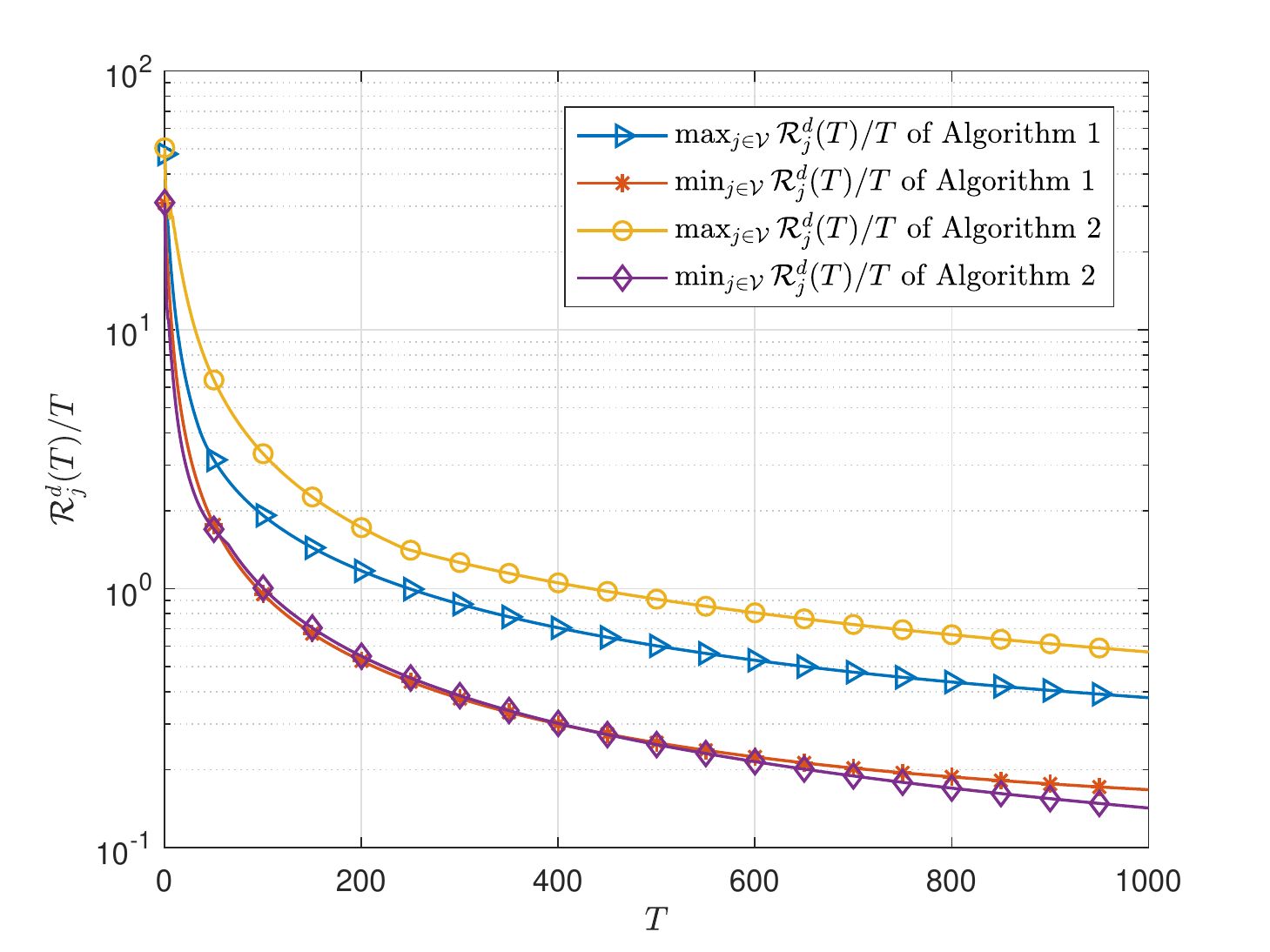} \vspace{-0ex}}
\caption{Comparisons of the maximum and minimum dynamic regret over time horizon $T$ 
between
Algorithm 1 and Algorithm 2.}
\label{figureexa1dd}
\end{figure}

Moreover, we compare the convergence performance between Algorithm 1 and Algorithm 2 under the same settings, except that Algorithm 1 adopts the local zeroth-order oracle  (\ref{GForacle})  with $\mu=10^{-3}$ and Algorithm 2 uses the explicit gradient of local cost functions.
The maximum and minimum average dynamic regrets over time horizon $T$, 
defined as 
i.e., $\max_{j\in\mathcal{V}}\mathcal{R}_{j}(T)/T$ and $\min_{j\in\mathcal{V}}\mathcal{R}_{j}(T)/T$, are depicted in Fig. \ref{figureexa1dd}.
The result is consistent with our theoretical results established in Theorem \ref{NetworkRegret} and Theorem \ref{Theo4G},
and further verifies the effectiveness of our zeroth-order algorithm
since the Algorithm 1 can achieve comparable performance with Algorithm 2 where explicit gradient information is used.

\subsection{Example II: Targets Tracking}
In this example, we consider a tracking problem where sensors collaborate to track moving targets. It is of practical significance to study distributed online algorithms over unbalanced networks  as communications between sensors are usually directed.
Moreover, dynamic regret is qualified to measure the tracking performance in contrast to static regret.

We consider a network consisting of $n=6$ nodes, and the information
sharing among the nodes is depicted in Fig. \ref{figurea}. 
All nodes aim to collaboratively track three time-varying signals via local communication.
Specifically, each signal $x^{\star}_{i,t}=[x^{\star}_{i1,t},x^{\star}_{i2,t}]^{\text{T}}$, $i\in\{1,2,3\}$ is described as follows 
\begin{eqnarray}\label{Sye}
x_{i,t}^{\star}=
\left[
  \begin{array}{ccc}
    x_{i1,t}^{\star}\\
    \dot{x}_{i1,t}^{\star}\\
  \end{array}
\right]
=
\left[
  \begin{array}{ccc}
    \kappa_{i}\sin(\omega_{i} t+\nu_{i})\\
   \omega_{i}\kappa_{i}\cos(\omega_{i} t+\nu_{i})\\
  \end{array}
\right]
\end{eqnarray}
where $x_{i1,t}^{\star}$ is the position of target $i$,  $\dot{x}_{i1,t}^{\star}$ denotes the velocity of target $i$ at time $t$, $\kappa_{i}$ is the amplitude, $\omega_{i}$ represents the angular frequency, and $\nu_{i}$ is the phase of the target $i$.
At time step $t$, each node $i\in\mathcal{V}$ observes $y_{i,t}$ via the measurement model $y_{i,t}=C_{i}x^{\star}_{t}$, where $x^{\star}_{t}=[(x^{\star}_{1,t})^{\text{T}},(x^{\star}_{2,t})^{\text{T}},(x^{\star}_{3,t})^{\text{T}}]^{\text{T}}$ and $C_{i}\in\mathbb{R}^{1\times 6}$ is the measurement matrix that is generated randomly.
To track the moving targets, each node communicates with its neighbors aiming to minimize the global cost function $f_{t}(x)=\frac{1}{2}\sum_{i=1}^{6}\Vert C_{i}x-y_{i,t}\Vert^{2}$.

\begin{figure}[ptbh]
\centerline{
\includegraphics[width=8cm]{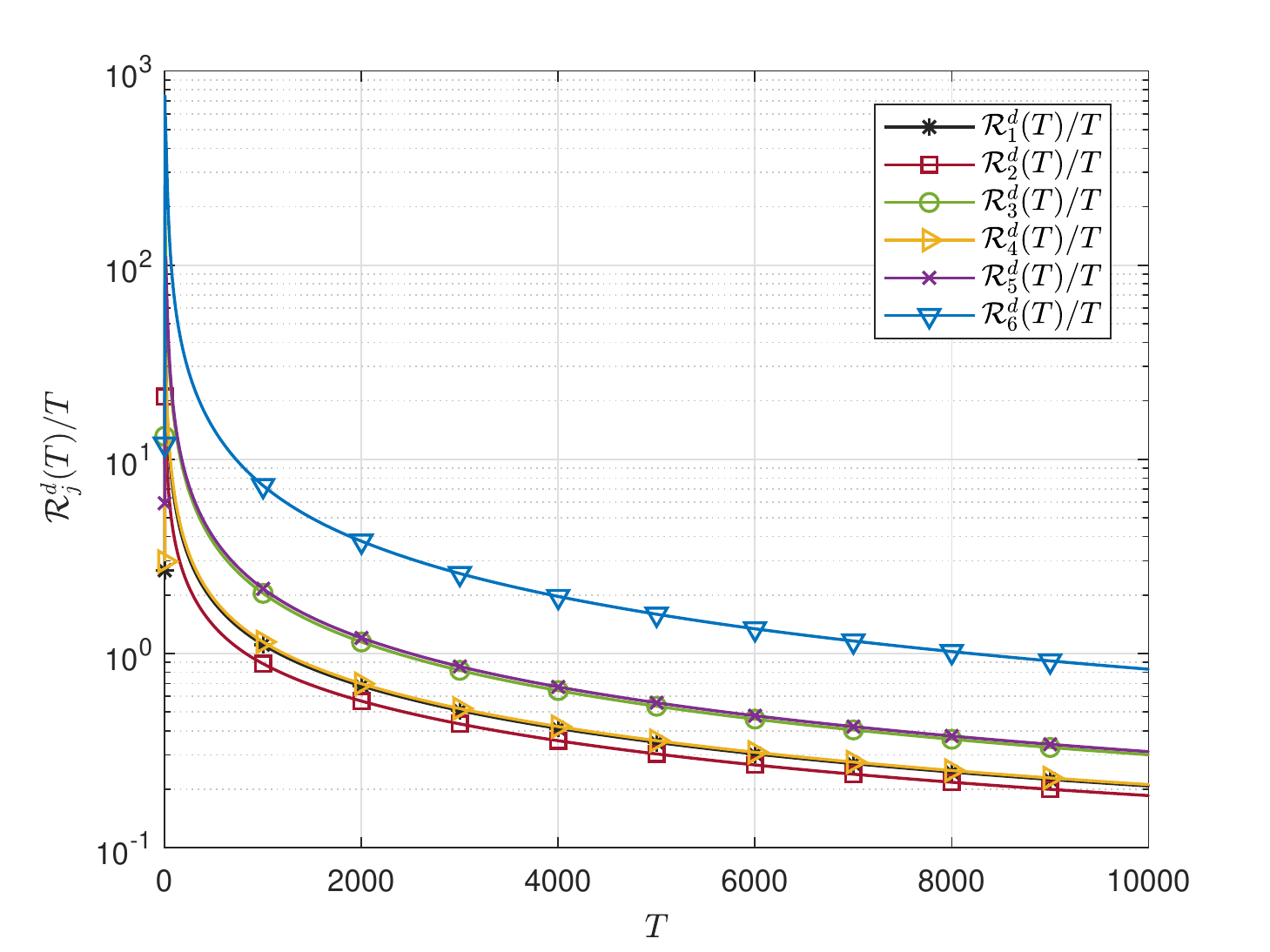} \vspace{-0ex}}
\caption{Trajectories of the average dynamic regret 
over time horizon $T$ 
of Algorithm 1.}
\label{figureAlg2Reg}
\end{figure}

\begin{figure}[ptbh]
\centerline{
\includegraphics[width=9cm]{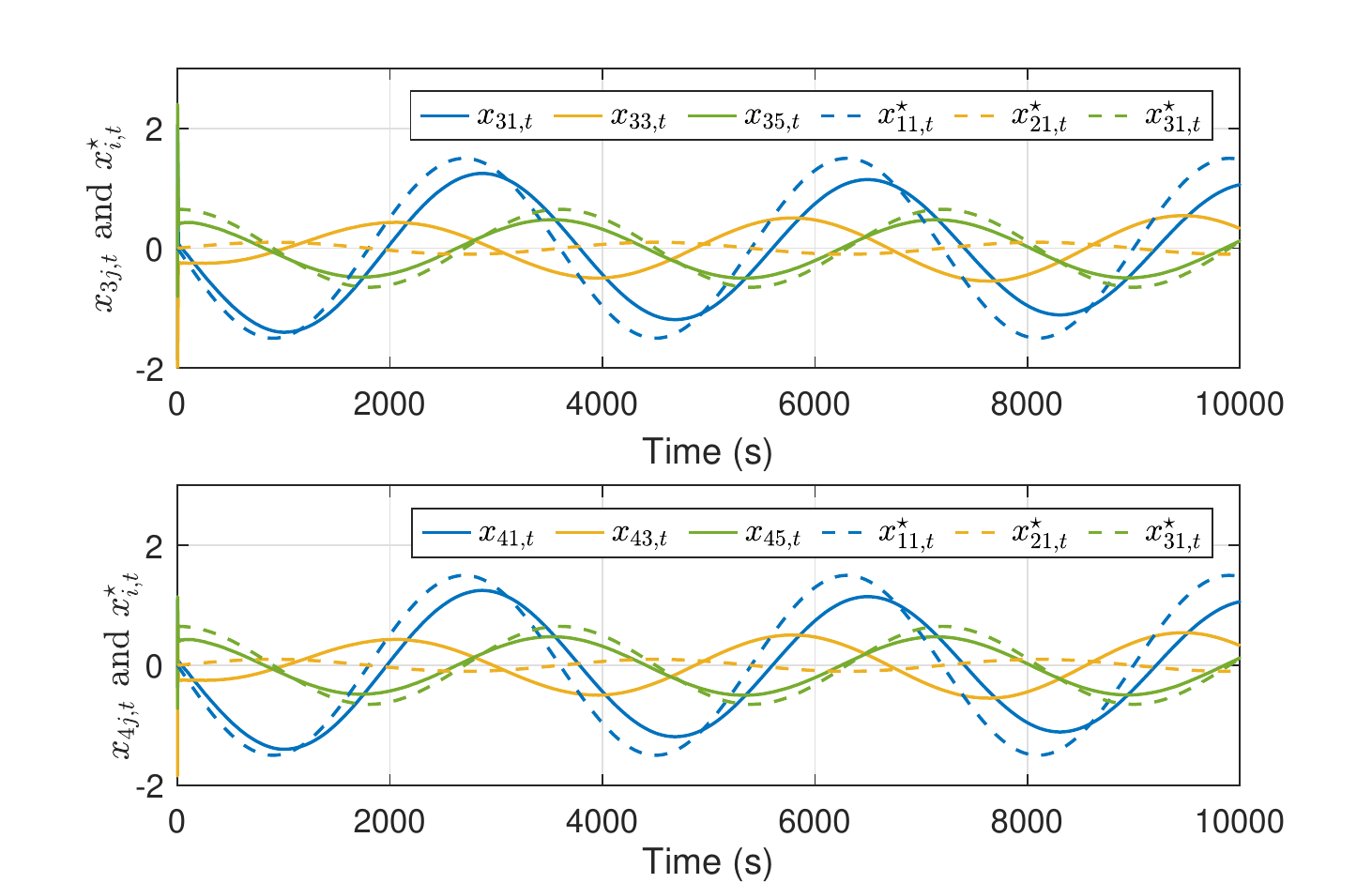} \vspace{-0ex}}
\caption{Trajectories of the targets $x_{i1,t}^{\star}$, $i\in\{1,2,3\}$, and the decision variables $x_{3j,t}$ and $x_{4j,t}$, $j\in\{1,3,5\}$, generated by Algorithm 1.}
\label{figureTrajAl1}
\end{figure}

In our simulations,  the amplitude $\kappa_{i}$ and the phase $\nu_{i}$, for $i\in\mathcal{V}$, are random variables drawn from uniform distribution in $[0,3]$ and $[0,\pi]$, respectively.
We let the sampling frequency be 100Hz.
Similar to the  previous example, we verify the performance of Algorithm 1 with $\alpha_{t}=\frac{1}{2\sqrt{t+1}}$ and $\mu=10^{-3}$. As shown in Fig. \ref{figureAlg2Reg}, the average  dynamic regret over $T$ calculated via any sequence of $\{x_{i,t}\}_{t=0}^{T}$, $i\in\mathcal{V}$, asymptotically converges 
at a sublinear rate,
which is consistent with our results  established in Theorem \ref{NetworkRegret}.
Furthermore, we show the trajectories of the targets $x_{i1,t}^{\star}$, $i\in\{1,2,3\}$, and the decisions $x_{3j,t}$ and $x_{4j,t}$, $j\in\{1,3,5\}$, in Fig. \ref{figureTrajAl1} with $x_{ij,t}$ being the $j$-th entry of the decision made by agent $i$ at $t$. It can be observed that node 3 and node 4 both can track the three moving targets within the targets' small neighborhood.
In fact, the tracking performance can be further improved by tuning the sampling frequency.
However, utilizing the diminishing step-size prevents our algorithms from tracking quickly moving targets. It is thus of interest to develop distributed online optimization algorithms in time-varying unbalanced networks that admits more aggressive step sizes to obtain better tracking performance, which remains to be considered in our future work.


\subsection{Example III: The Dynamic Sparse Recovery Problem}
\begin{figure}[ptbh]
\centerline{
\includegraphics[width=8cm]{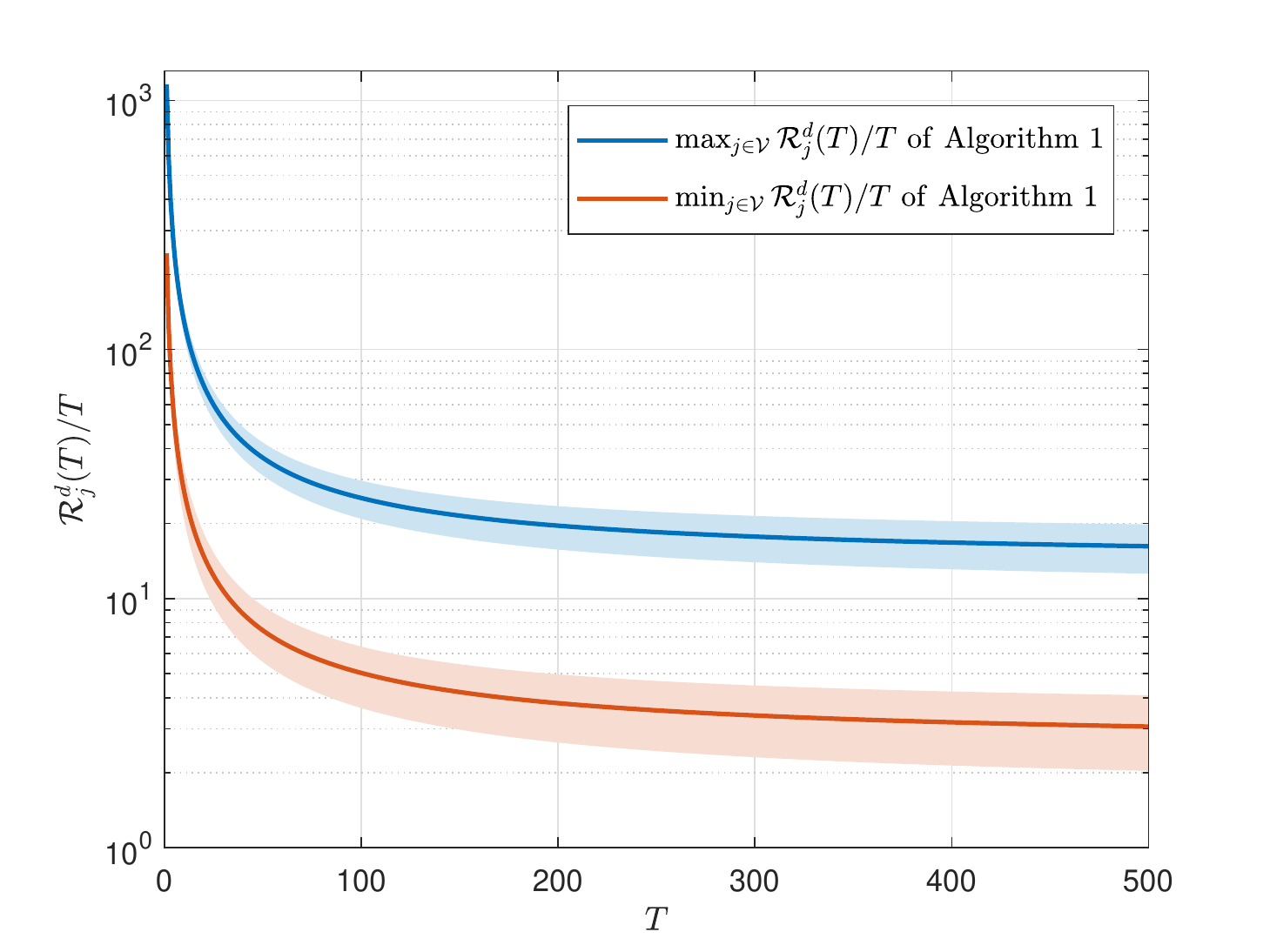} \vspace{-0ex}}
\caption{{\color{black}Trajectories of the maximum and minimum dynamic regret over time horizon $T$ of Algorithm 1.}}
\label{figureexa3_1}
\end{figure}

In this example, we consider the dynamic sparse signal recovery problem with the goal of estimating a time-varying sparse parameter in a distributed manner. This problem has been widely investigated in signal processing literature. 
Specifically, we compare our algorithm with the ones in \cite{ShahrampourS2018Distributed,ZhangY2019A,Luk2019Online}. 
Note that all of the algorithms in these works can be only applicable to (fixed or time-varying) balanced networks. 
%
To investigate these algorithms within a unified framework, we consider a network with $n$ sensors, which can be represented by a complete graph. All sensors collectively solve the following optimization problem: 
\begin{eqnarray}
x_{t}^{\star}=\arg\min_{x} f_{t}(x)
\end{eqnarray}
with  
\begin{eqnarray}
f_{t}(x)&\hspace*{-0.5em}=\hspace*{-0.5em}&\frac{1}{n}\sum_{i=1}^{n}\left(\left\Vert z_{i,t}-C_{i,t}x\right\Vert_{2}^{2}+\gamma\Vert x\Vert_{2}^{2}+\sigma\Vert x\Vert_{1}\right),  
\end{eqnarray}
where $x\in\mathbb{R}^{m}$, $C_{i,t}\in\mathbb{R}^{d\times m}$ is the observation matrix of sensor $i$ at round $t$, 
$\gamma$ and $\sigma$ are constant regularization parameters used to prevent overfitting, 
and $z_{i,t}\in\mathbb{R}^{d}$ is the measurement given by 
\begin{eqnarray}
z_{i,t}&\hspace*{-0.5em}=\hspace*{-0.5em}&C_{i,t}\varpi_{t}+\vartheta_{i,t}
\end{eqnarray}
with 
$\varpi_{t}\in\mathbb{R}^{m}$ being the time-varying sparse signal of interest and 
$\vartheta_{i,t}\in\mathbb{R}^{d}$ being the noise.  
In this example, 
we set $n=40$, $d=3$, $m=8$, $\gamma=\frac{1}{100d^{2}n}$ and $\sigma=\frac{1}{20d}$. 
The initial state of $\varpi_{t}$ is chosen to be a sparse vector with $2$ entries of value $1$, and all other entries equal to zero. 
Let $\mathcal{S}_{t}=\{\kappa|[\varpi_{t}]_{\kappa}>0\}$ be the support of $\varpi_{t}$, which is updated as follows: 
\begin{equation} \label{Support}
\mathcal{S}_{t+1}=
\begin{cases}
\mathcal{S}_{t},  &\text{with probability} ~1-1/t \\
\{\mathcal{S}_{t}\backslash \{\tau_{t}\}\}\cup\{\tau_{t}'\}, &\text{with probability}~ 1/t
\end{cases}
\end{equation}
where $\tau_{t}$ and $\tau_{t}'$ are randomly chosen from $\mathcal{S}_{t}$ and the set $\{1,2,...,m\}\backslash\mathcal{S}_{t}$, respectively. 
We add noise $\varsigma_{t}$ to $\varpi_{t}$ and then normalize the obtained vector so that 
\begin{eqnarray}
\varpi_{t+1}&\hspace*{-0.5em}=\hspace*{-0.5em}&\frac{\varpi_{t}+\varsigma_{t}}{\Vert \varpi_{t}+\varsigma_{t}\Vert}
\end{eqnarray}
where $[\varsigma_{t}]_{\kappa}\sim \mathcal{N}(0,1/t^{2})$ for $\kappa\in\mathcal{S}_{t+1}$ and $[\varsigma_{t}]_{\kappa}=0$ otherwise. 
By this means, the non-zero entries of $\varpi_{t}$ are time-varying and their variations decay over time. 
The related parameters are properly selected such that the minimizer sequence $\{x_{t}^{\star}\}_{t=1}^{T}$ varies slowly.  
We verify the performance of Algorithm 1 with $\alpha_{t}=\frac{1}{\sqrt{t+1}}$ and $\mu=10^{-3}$. {\color{black}We run Algorithm 1 for 100 trials. As shown in Fig. \ref{figureexa3_1}, the average maximum and minimum dynamic regrets over $T$ both converge sublinearly, which are consistent with our theoretical results.} 

\begin{figure}[ptbh]
\centerline{
\includegraphics[width=8cm]{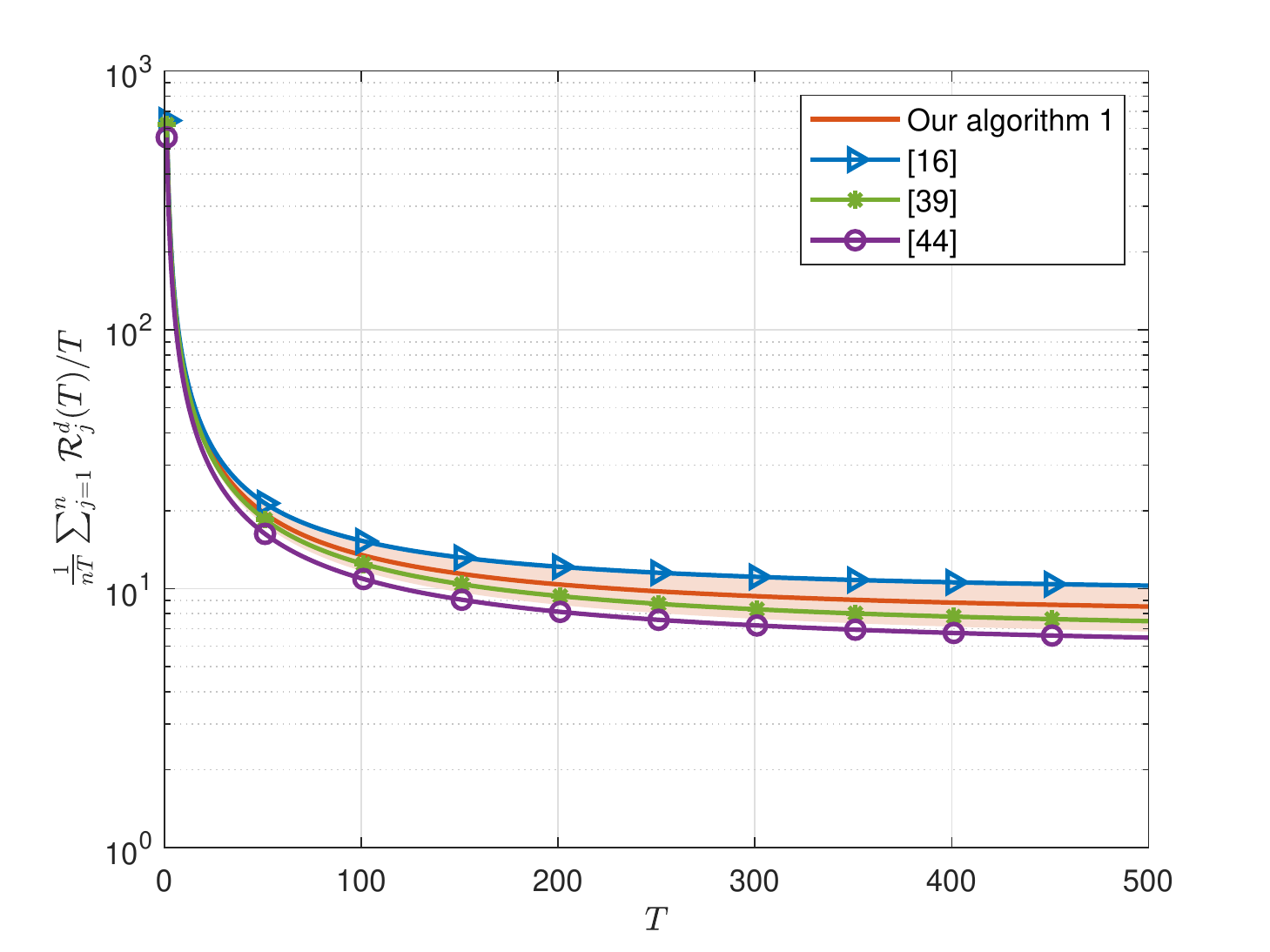} \vspace{-0ex}}
\caption{{\color{black}Comparison of the time-average dynamic regret averaged over all nodes between our algorithms and the existing ones.}}
\label{figureexa3_algorithms}
\end{figure}

We further compare our algorithm against the ones in \cite{ShahrampourS2018Distributed,ZhangY2019A,Luk2019Online}. 
For fair comparison, we adopt the 
average of $\frac{\mathcal{R}_{j}^{d}(T)}{T}$ over all nodes $j\in\mathcal{V}$
as the performance metric of the whole network. 
The results are depicted in Fig. \ref{figureexa3_algorithms}, which reflect that 
the proposed Algorithm 1 
can achieve comparable performance with the existing algorithms despite the presence of gradient estimate errors. 
Note that our algorithm can be applied to a broader family of networks since it does not require networks to be fixed or balanced in contrast to \cite{ShahrampourS2018Distributed,ZhangY2019A,Luk2019Online}. 

\begin{figure}[ptbh]
\centerline{
\includegraphics[width=8cm]{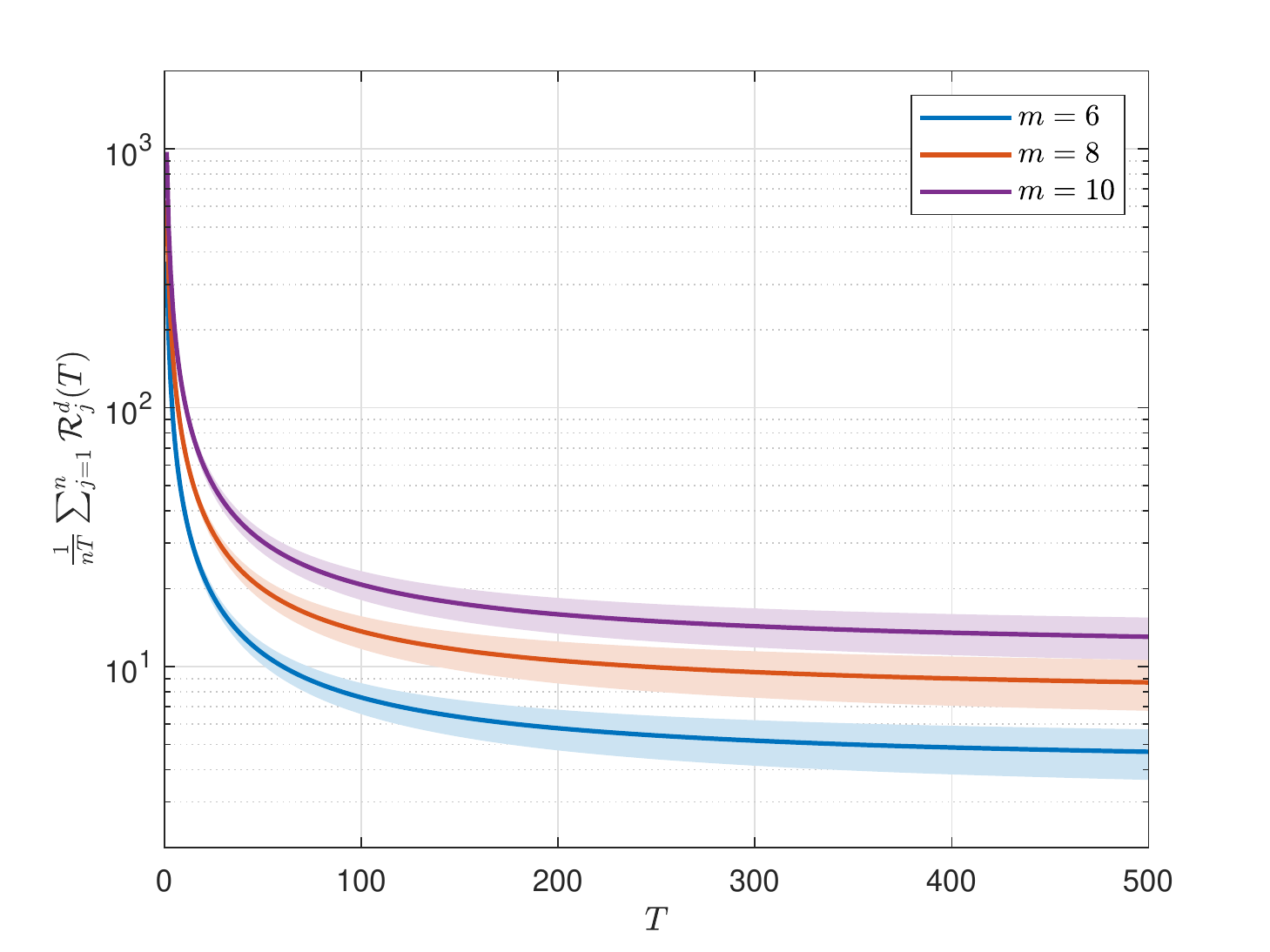} \vspace{-0ex}}
\caption{{\color{black}The influence of the problem dimension $m$ on the performance of Algorithm 1. 
}}
\label{figureexa3_ddimens}
\end{figure}

Finally, considering the fact that zeroth-order optimization algorithms usually suffer from deterioration in performance as the problem dimension increases, we investigate the influence of the problem dimension on the average dynamic regret of our algorithm. 
We run our algorithm for three different choices of the problem dimension $m$, i.e., $m=6$, $m=8$ and $m=10$, respectively.  
The results are depicted in Fig. \ref{figureexa3_ddimens}. It can be seen that the proposed zeroth-order optimization algorithm indeed achieves better performance with smaller problem dimension.

%
%
%
%
%
%
%
%
%

\section{Conclusions}\label{sec-VI}

In this paper, we have investigated the distributed online constrained optimization problems over time-varying unbalanced
digraphs without explicit subgradients. To cope with the problem, a novel consensus-based distributed online optimization algorithm with a local randomized zeroth-order oracle has been proposed. The dynamic regret bound of the algorithm has been characterized under mild conditions, which showed that the algorithm
can achieve comparable performance with its counterpart subgradient-based algorithm.
Finally, simulations on distributed target tracking problem and dynamic sparse signal recovery problem in sensor networks have been conducted to verify the effectiveness of proposed algorithms.
Future works can focus on
deriving a tighter dynamic regret bound,
and considering the scenario where inequalities constraints are involved in the problem.
It is also of interest to 
 develop dimension-insensitive zeroth-order distributed online optimization algorithms under mild assumptions. 

\section*{Appendix A}

\textit{Proof of Lemma \ref{Rescale}.}
Recalling $\phi_{i,t+1}=\sum_{j=1}^{n}a_{ij,t}\phi_{j,t}$ and $b_{ij,t}=\frac{1}{\phi_{i,t+1}}a_{ij,t}\phi_{j,t}$, along with the fact that $\mathcal{A}_{t}$ is a column matrix for $ t\geq 0$, it is not hard to verify that $\phi_{t+1}^{\text{T}}\mathcal{B}_{t}=\phi_{t}^{\text{T}}$.
Then, conducting the mathematical induction for $\phi_{t}$ yields that 
\begin{eqnarray}
\phi_{t+1}^{\text{T}}=\phi_{T}^{\text{T}}\mathcal{B}(T:t+1)
\end{eqnarray}
On the other hand, by exploiting the column stochasticity of $\mathcal{A}_{t}$, we have 
\begin{eqnarray}
\bm{1}_{n}^{\text{T}}\phi_{t+1}=\bm{1}_{n}^{\text{T}}\mathcal{A}_{t}\phi_{t}=...=\sum_{i=1}^{n}\phi_{i,0}=n
\end{eqnarray}
%
%
%
Therefore, we can further obtain that
\begin{eqnarray}
&&\left|\frac{1}{n}\phi_{i,t+1}-\pi_{i,t+1}\right|    \notag\\
&&=
\left| \frac{1}{n}\sum_{l=1}^{n}\phi_{l,T}b_{li}(T:t+1)-\frac{\pi_{i,t+1}}{n}\sum_{l=1}^{n}\phi_{l,T} \right|  \notag\\
&&\leq
\frac{1}{n}\sum_{l=1}^{n}\phi_{l,T}|b_{li}(T:t+1)-\pi_{i,t+1}| \notag\\
&&\leq
C\lambda^{T-t-1}
\end{eqnarray}
where
Lemma \ref{Convergence-z}(a) has been utilized to obtain the last inequality.
The proof is thus completed.
$\hfill \blacksquare$

\section*{Appendix B}

\textit{Proof of Lemma \ref{theoremcons}.}
Define 
\begin{eqnarray}
v_{i,t}&\hspace*{-0.5em}\triangleq \hspace*{-0.5em}&\sum_{j=1}^{n}b_{ij,t}x_{j,t}  \notag\\
w_{i,t}&\hspace*{-0.5em}\triangleq \hspace*{-0.5em}&\mathcal{P}_{{\color{black}(1-\xi)\Omega}}\left(v_{i,t}-\frac{1}{\phi_{i,t+1}}\alpha_{t}{\color{black}\hat{g}_{i,t}(x_{i,t})}\right)-v_{i,t}. \notag
\end{eqnarray}
%
%
Then (\ref{Alg_2}) can be rewritten as the following perturbed form:
\begin{eqnarray} \label{Individual-1}
x_{i,t+1}
&\hspace*{-0.5em}=\hspace*{-0.5em}&
\sum_{j=1}^{n}b_{ij,t}x_{j,t}+w_{i,t}.
\end{eqnarray}
Conducting the mathematical induction for (\ref{Individual-1}) yields that
\begin{eqnarray} \label{Total-1}
x_{i,t+1}
&\hspace*{-0.5em}=\hspace*{-0.5em}&
\sum_{j=1}^{n}[\mathcal{B}(t+1:0)]_{ij}x_{j,0}   \notag\\
&&+\sum_{l=1}^{t+1}\sum_{j=1}^{n}[\mathcal{B}(t+1:l)]_{ij}w_{j,l-1}.
\end{eqnarray}

On the other hand, multiplying both sides of (\ref{Individual-1}) by $\pi_{i,t+1}$ and then summing the obtained equality over $i\in\mathcal{V}$ yields that 
\begin{eqnarray} \label{Bar-x}
\bar{x}_{t+1}
&\hspace*{-0.5em}=\hspace*{-0.5em}&
\bar{x}_{t}+\sum_{i=1}^{n}\pi_{i,t+1}w_{i,t},
\end{eqnarray}
where Lemma \ref{Convergence-z}(c) has been exploited to obtain this equality.
Performing summations on (\ref{Bar-x}) over $t\in\mathbb{N}$ leads to
\begin{eqnarray} \label{TotalBar_x}
\bar{x}_{t+1}
&\hspace*{-0.5em}=\hspace*{-0.5em}&
\bar{x}_{0}+\sum_{l=1}^{t+1}\sum_{i=1}^{n}\pi_{i,l}w_{i,l-1}.
\end{eqnarray}
Then, combining (\ref{Total-1}) and (\ref{TotalBar_x}) along with Lemma \ref{Convergence-z}(a), it gives that
\begin{eqnarray} \label{Difference1}
&&\hspace{-1.5em}\Vert x_{i,t+1}-\bar{x}_{t+1}\Vert   \notag\\
&&\hspace{-1.5em}=
\sum_{j=1}^{n}\left|[\mathcal{B}(t+1:0)]_{ij}-\pi_{j,0}\right|\Vert x_{j,0}\Vert  \notag\\
&&\hspace{-1.5em}~~~+\sum_{l=1}^{t+1}\sum_{j=1}^{n}\left|[\mathcal{B}(t+1:l)]_{ij}-\pi_{j,l}\right|\Vert w_{j,l-1}\Vert \notag\\
&&\hspace{-1.5em}\leq
C\lambda^{t+1}\sum_{j=1}^{n}\Vert x_{j,0}\Vert +C\sum_{l=1}^{t+1}\lambda^{t+1-l}\sum_{j=1}^{n}\Vert w_{j,l-1} \Vert
\end{eqnarray}
Now, it remains to bound $\Vert w_{i,t}\Vert$.
Note that $v_{i,t}\in\Omega$ is a convex combination of $x_{i,t}\in\Omega$, it thus follows that
\begin{eqnarray} \label{wit}
\Vert w_{i,t}\Vert
&\hspace*{-0.5em}\leq\hspace*{-0.5em}&
\left\Vert v_{i,t}-\frac{\alpha_{t}{\color{black}\hat{g}_{i,t}(x_{i,t})}}{\phi_{i,t+1}}-v_{i,t}\right\Vert \notag\\
&\hspace*{-0.5em}\leq\hspace*{-0.5em}&
\theta\alpha_{t}\Vert {\color{black}\hat{g}_{i,t}(x_{i,t})} \Vert \notag\\
&\hspace*{-0.5em}\leq\hspace*{-0.5em}&
m\theta G\alpha_{t} 
\end{eqnarray}
where the non-expansive projection property (see e.g., \cite{nedic2010Constrained}) has been exploited to obtain the first inequality, 
and the last two inequalities follows from Lemma \ref{BoundLe} and Lemma \ref{gradient-freePro}, respectively. 
It is straightforward to obtain (\ref{Theorem-1}) by combining (\ref{Difference1}) and (\ref{wit}). 
%
%
%
%
$\hfill \blacksquare$


\section*{Appendix C}

\textit{Proof of Theorem \ref{Lemma_4}.}
{\color{black}Define $\tilde{x}_{t}^{\star}\triangleq (1-\xi)x_{t}^{\star}$. Obviously, $\tilde{x}_{t}^{\star}\in(1-\xi)\Omega$, 
and $\bar{x}_{t}\in(1-\xi)\Omega$. 
We further define an auxiliary function $D_{t}$ as 
\begin{eqnarray}
D_{t}\triangleq \frac{1}{2}\left(\bar{x}_{t}-\tilde{x}_{t}^{\star}\right)^{\text{T}}\left(\bar{x}_{t}-\tilde{x}_{t}^{\star}\right)
\end{eqnarray} 
Where $\bar{x}_{t}$ has been defined in (\ref{x_bar-1}).}
Then, we consider the variation of $D_{t}$ as follows
\begin{eqnarray} \label{variation}
\Delta D_{t}
&\hspace*{-0.5em}\triangleq \hspace*{-0.5em}&
D_{t+1}-D_{t}\notag\\
&\hspace*{-0.5em}=\hspace*{-0.5em}&
-\frac{1}{2}\Vert \bar{x}_{t+1}-\bar{x}_{t}\Vert^{2}+\frac{1}{2}\langle {\color{black}\tilde{x}_{t+1}^{\star}}+{\color{black}\tilde{x}_{t}^{\star}}-2\bar{x}_{t+1},  \notag\\
&&{\color{black}\tilde{x}_{t+1}^{\star}}-{\color{black}\tilde{x}_{t}^{\star}}\rangle
+\langle \bar{x}_{t+1}-\bar{x}_{t},\bar{x}_{t+1}-{\color{black}\tilde{x}_{t}^{\star}}\rangle   \notag\\
&\hspace*{-0.5em}\leq\hspace*{-0.5em}&
-\frac{1}{2}\Vert \bar{x}_{t+1}-\bar{x}_{t}\Vert^{2}
+2R\Vert {\color{black}\tilde{x}_{t+1}^{\star}}-{\color{black}\tilde{x}_{t}^{\star}}\Vert      \notag\\
&&+\langle \bar{x}_{t+1}-\bar{x}_{t},\bar{x}_{t+1}-{\color{black}\tilde{x}_{t}^{\star}}\rangle
\end{eqnarray}
where the inequality follows from the boundedness of the constraint set.
According to (\ref{Bar-x}), the last term in (\ref{variation}) can be further expanded as follows
\begin{eqnarray} \label{angle}
&&\langle \bar{x}_{t+1}-\bar{x}_{t},\bar{x}_{t+1}-{\color{black}\tilde{x}_{t}^{\star}}\rangle \notag\\
&&=
\left\langle \sum_{i=1}^{n}\pi_{i,t+1}w_{i,t},\bar{x}_{t+1}-{\color{black}\tilde{x}_{t}^{\star}}\right\rangle   \notag\\
&&=
\sum_{i=1}^{n}\pi_{i,t+1}\left\langle w_{i,t},\bar{x}_{t+1}-x_{i,t+1}\right\rangle  \notag\\
&&~~~+\sum_{i=1}^{n}\pi_{i,t+1}\left\langle w_{i,t},x_{i,t+1}-{\color{black}\tilde{x}_{t}^{\star}}\right\rangle
\end{eqnarray}

For the first term on the right side of (\ref{angle}), we have
\begin{eqnarray}\label{Firstterm}
&& \sum_{i=1}^{n}\pi_{i,t+1}\left\langle w_{i,t},\bar{x}_{t+1}-x_{i,t+1}\right\rangle  \notag\\
&&\leq
\sum_{i=1}^{n}\pi_{i,t+1}\Vert w_{i,t}\Vert \Vert\bar{x}_{t+1}-x_{i,t+1}\Vert
\end{eqnarray}
Recalling (\ref{Difference1})-(\ref{wit}), 
and taking expectation on both sides of (\ref{Firstterm}),
we can obtain
\begin{eqnarray}\label{FirstEXP}
&& \mathbb{E}\left[\sum_{i=1}^{n}\pi_{i,t+1}\left\langle w_{i,t},\bar{x}_{t+1}-x_{i,t+1}\right\rangle \right]  \notag\\
&&\leq
{\color{black}m}\theta CG\alpha_{t}\lambda^{t+1}\sum_{j=1}^{n}\Vert x_{j,0}\Vert  \notag\\
&&~~~+n{\color{black}m^{2}}\theta^{2}CG^{2}\alpha_{t}\sum_{l=1}^{t+1}\lambda^{t+1-l}\alpha_{l-1}
\end{eqnarray}
%
where 
the fact that $\sum_{i=1}^{n}\pi_{i,t+1}=1$ has also been exploited to obtain the above inequality.

For the last term on the right side of (\ref{angle}), we have
\begin{eqnarray}\label{cross1}
&& \sum_{i=1}^{n}\pi_{i,t+1}\left\langle w_{i,t},x_{i,t+1}-{\color{black}\tilde{x}_{t}^{\star}}\right\rangle     \notag\\
&&=
 \sum_{i=1}^{n}\pi_{i,t+1}\left\langle w_{i,t}+\frac{\alpha_{t}{\color{black}\hat{g}_{i,t}(x_{i,t})}}{\phi_{i,t+1}},x_{i,t+1}-{\color{black}\tilde{x}^{\star}_{t}}\right\rangle   \notag\\
&&~~~ +\sum_{i=1}^{n}\pi_{i,t+1}\left\langle \frac{\alpha_{t}{\color{black}\hat{g}_{i,t}(x_{i,t})}}{\phi_{i,t+1}}, {\color{black}\tilde{x}^{\star}_{t}}-x_{i,t+1}\right\rangle
\end{eqnarray}
From Lemma 1 in \cite{nedic2010Constrained}, we know $(\mathcal{P}_{{\color{black}(1-\xi)\Omega}}(x)-x)^{\text{T}}(\mathcal{P}_{{\color{black}(1-\xi)\Omega}}(x)-y)\leq 0$, for all $x\in\mathbb{R}^{m}$ and $y\in{\color{black}(1-\xi)\Omega}$. Hence,
it can be observed that
\begin{eqnarray}\label{ProjPro}
\left\langle w_{i,t}+\frac{\alpha_{t}{\color{black}\hat{g}_{i,t}(x_{i,t})}}{\phi_{i,t+1}},x_{i,t+1}-{\color{black}\tilde{x}^{\star}_{t}}\right\rangle\leq 0
\end{eqnarray}
That is to say, the first term on the right side of (\ref{cross1}) is nonpositive.
Now, we turn to the last term on the right side of (\ref{cross1}), which can be calculated as follows
\begin{eqnarray}\label{cross2}
&&\sum_{i=1}^{n}\pi_{i,t+1}\left\langle \frac{\alpha_{t}{\color{black}\hat{g}_{i,t}(x_{i,t})}}{\phi_{i,t+1}}, {\color{black}\tilde{x}^{\star}_{t}}-x_{i,t+1}\right\rangle  \notag\\
&&=
\sum_{i=1}^{n}\pi_{i,t+1}\left\langle \frac{\alpha_{t}{\color{black}\hat{g}_{i,t}(x_{i,t})}}{\phi_{i,t+1}}, {\color{black}\tilde{x}^{\star}_{t}}-\bar{x}_{t} \right\rangle     \notag\\
&&~~~+\sum_{i=1}^{n}\pi_{i,t+1}\left\langle \frac{\alpha_{t}{\color{black}\hat{g}_{i,t}(x_{i,t})}}{\phi_{i,t+1}}, \bar{x}_{t}-\bar{x}_{t+1} \right\rangle   \notag\\
&&~~~+\sum_{i=1}^{n}\pi_{i,t+1}\left\langle \frac{\alpha_{t}{\color{black}\hat{g}_{i,t}(x_{i,t})}}{\phi_{i,t+1}}, \bar{x}_{t+1}-x_{i,t+1} \right\rangle
\end{eqnarray}
For the first term on the right side of (\ref{cross2}), it can be calculated as follows
\begin{eqnarray}\label{cross3}
&&\sum_{i=1}^{n}\pi_{i,t+1}\left\langle \frac{\alpha_{t}{\color{black}\hat{g}_{i,t}(x_{i,t})}}{\phi_{i,t+1}}, {\color{black}\tilde{x}^{\star}_{t}}-\bar{x}_{t} \right\rangle     \notag\\
&&=
\sum_{i=1}^{n}\pi_{i,t+1}\alpha_{t} \left\langle \frac{{\color{black}\hat{g}_{i,t}(x_{i,t})}-{\color{black}\nabla \hat{f}_{i,t}(\bar{x}_{t})}
}{\phi_{i,t+1}}, {\color{black}\tilde{x}^{\star}_{t}}-\bar{x}_{t} \right\rangle  \notag\\
&&~~~+\sum_{i=1}^{n}\pi_{i,t+1}\alpha_{t} \left\langle \frac{{\color{black}\nabla \hat{f}_{i,t}(\bar{x}_{t})}}{\phi_{i,t+1}}, {\color{black}\tilde{x}^{\star}_{t}}-\bar{x}_{t} \right\rangle
\end{eqnarray}
Taking expectation on the above equality with respect to $\mathcal{F}_{t-1}$,
we can obtain
\begin{eqnarray}\label{EXPcros}
&&\mathbb{E}\left[  \sum_{i=1}^{n}\pi_{i,t+1}\left\langle \frac{\alpha_{t}{\color{black}\hat{g}_{i,t}(x_{i,t})}}{\phi_{i,t+1}}, {\color{black}\tilde{x}^{\star}_{t}}-\bar{x}_{t}\right\rangle \big | \mathcal{F}_{t-1} \right]   \notag\\
&&=
\sum_{i=1}^{n}\pi_{i,t+1}\alpha_{t} \left\langle \frac{
{\color{black}\nabla \hat{f}_{i,t}(x_{i,t})-\nabla \hat{f}_{i,t}(\bar{x}_{t})}
}{\phi_{i,t+1}}, {\color{black}\tilde{x}^{\star}_{t}}-\bar{x}_{t} \right\rangle  \notag\\
&&~~~+\sum_{i=1}^{n}\pi_{i,t+1}\alpha_{t} \left\langle \frac{{\color{black}\nabla \hat{f}_{i,t}(\bar{x}_{t})}}{\phi_{i,t+1}}, {\color{black}\tilde{x}^{\star}_{t}}-\bar{x}_{t} \right\rangle  \notag\\
&&\leq
\theta \alpha_{t}\sum_{i=1}^{n}\pi_{i,t+1} 
\left\Vert {\color{black}\nabla \hat{f}_{i,t}(x_{i,t})-\nabla \hat{f}_{i,t}(\bar{x}_{t})}\right\Vert \Vert {\color{black}\tilde{x}^{\star}_{t}}-\bar{x}_{t} \Vert  \notag\\
&&~~~+\sum_{i=1}^{n}\pi_{i,t+1}\alpha_{t} \left\langle \frac{{\color{black}\nabla \hat{f}_{i,t}(\bar{x}_{t})}}{\phi_{i,t+1}}, {\color{black}\tilde{x}^{\star}_{t}}-\bar{x}_{t} \right\rangle
\end{eqnarray}
where Lemma \ref{gradient-freePro}(b) has been employed to obtain the equality, while the inequality follows from Lemma \ref{BoundLe}. Note that $\Vert {\color{black}\tilde{x}^{\star}_{t}}-\bar{x}_{t} \Vert\leq 2R$. In addition, Lemma \ref{gradient-freePro} implies that $ {\color{black}\Vert \nabla \hat{f}_{i,t}(x_{i,t})-\nabla \hat{f}_{i,t}(\bar{x}_{t})\Vert\leq L\Vert x_{i,t} - \bar{x}_{t}\Vert}$. Thus, (\ref{EXPcros}) can be further manipulated as
\begin{eqnarray}\label{EXPcrosImp}
&&\mathbb{E}\left[  \sum_{i=1}^{n}\pi_{i,t+1}\left\langle \frac{\alpha_{t}{\color{black}\hat{g}_{i,t}(x_{i,t})}}{\phi_{i,t+1}}, {\color{black}\tilde{x}^{\star}_{t}}-\bar{x}_{t}\right\rangle \big | \mathcal{F}_{t-1} \right]   \notag\\
&&\leq
2\theta RL\alpha_{t}\sum_{i=1}^{n}\pi_{i,t+1}\Vert x_{i,t}-\bar{x}_{t}\Vert \notag\\
&&~~~+\alpha_{t}\sum_{i=1}^{n}\frac{\pi_{i,t+1}}{\phi_{i,t+1}}\left({\color{black}\hat{f}_{i,t}(\tilde{x}_{t}^{\star})}
-f_{i,t}(\bar{x}_{t})\right)
\end{eqnarray}
where the convexity of {\color{black}$\hat{f}_{i,t}(x)$} and the fact that $f_{i,t}(x)\leq {\color{black}\hat{f}_{i,t}(x)}$ for all $i\in\mathcal{V}$ and $t\geq 0$ in Lemma \ref{gradient-freePro} (a) have been exploited to obtain the inequality.
By recalling (\ref{Bar-x}), the second term on the right side of (\ref{cross2}) can be bounded as
\begin{eqnarray}\label{cross4}
&&\hspace{-1.5em}\sum_{i=1}^{n}\pi_{i,t+1}\left\langle \frac{\alpha_{t}{\color{black}\hat{g}_{i,t}(x_{i,t})}}{\phi_{i,t+1}}, \bar{x}_{t}-\bar{x}_{t+1} \right\rangle     \notag\\
&&\hspace{-1.5em}\leq
\frac{1}{2}\left\Vert \sum_{i=1}^{n}\pi_{i,t+1}\frac{\alpha_{t}{\color{black}\hat{g}_{i,t}(x_{i,t})}}{\phi_{i,t+1}}\right\Vert^{2}
+\frac{1}{2}\left\Vert \bar{x}_{t}- \bar{x}_{t+1}\right\Vert^{2}   \notag\\
&&\hspace{-1.5em}\leq
\frac{1}{2}\theta^{2}\alpha_{t}^{2}\sum_{i=1}^{n}\pi_{i,t+1}\Vert {\color{black}\hat{g}_{i,t}(x_{i,t})}\Vert^{2}
+\frac{1}{2}\left\Vert \bar{x}_{t}- \bar{x}_{t+1}\right\Vert^{2}
\end{eqnarray}
where the last inequality follows from 
Lemma \ref{BoundLe}, 
along with the fact that 
$\Vert \sum_{i=1}^{n}\pi_{i,t}x_{i}\Vert^{2}\leq\sum_{i=1}^{n}\pi_{i,t}\Vert x_{i}\Vert^{2}$ for all $x_{i}\in\mathbb{R}^{m}$.  
Taking expectation on both sides of (\ref{cross4}) 
yields that
\begin{eqnarray}\label{EXPcross4}
&&\mathbb{E}\left[ \sum_{i=1}^{n}\pi_{i,t+1}\left\langle \frac{\alpha_{t}{\color{black}\hat{g}_{i,t}(x_{i,t})}}{\phi_{i,t+1}}, \bar{x}_{t}-\bar{x}_{t+1} \right\rangle \right]  \notag\\
&&\leq
\frac{1}{2}{\color{black}m}^{2}\theta^{2}G^{2}\alpha_{t}^{2}+\frac{1}{2}\mathbb{E}\left[ \Vert \bar{x}_{t}
- \bar{x}_{t+1}\Vert^{2}\right]
\end{eqnarray}
The third term on the right side of (\ref{cross2}) can be bounded as
\begin{eqnarray}\label{cross5}
&&\sum_{i=1}^{n}\pi_{i,t+1}\left\langle \frac{\alpha_{t}{\color{black}\hat{g}_{i,t}(x_{i,t})}}{\phi_{i,t+1}}, \bar{x}_{t+1}-x_{i,t+1} \right\rangle     \notag\\
&&\leq
\theta \alpha_{t}\sum_{i=1}^{n}\pi_{i,t+1}\Vert {\color{black}\hat{g}_{i,t}(x_{i,t})}\Vert\Vert x_{i,t+1}-\bar{x}_{t+1}\Vert
\end{eqnarray}
Recalling (\ref{Difference1})-(\ref{wit}), and taking expectation on both sides of (\ref{cross5}), we obtain 
\begin{eqnarray}\label{Expcross5}
&& \mathbb{E}\left[\sum_{i=1}^{n}\pi_{i,t+1}\left\langle \frac{\alpha_{t}{\color{black}\hat{g}_{i,t}(x_{i,t})}}{\phi_{i,t+1}}, \bar{x}_{t+1}-x_{i,t+1} \right\rangle\right]  \notag\\
&&\leq
n{\color{black}m}^{2}\theta^{2}CG^{2}\alpha_{t}\sum_{l=1}^{t+1}\lambda^{t+1-l}\alpha_{l-1}  \notag\\
&&~~~+{\color{black}m}\theta  CG\alpha_{t}\lambda^{t+1}\sum_{j=1}^{n}\Vert x_{j,0}\Vert. 
\end{eqnarray}

Now, taking expectation on both sides of (\ref{variation}) and applying the above relations, we can obtain
%
\begin{eqnarray}\label{Expvarisu}
&&\mathbb{E}[\Delta D_{t}] \notag\\
&&\leq
\alpha_{t}\sum_{i=1}^{n}\frac{\pi_{i,t+1}}{\phi_{i,t+1}}\mathbb{E}\left[{\color{black}\hat{f}_{i,t}(\tilde{x}_{t}^{\star})}
-f_{i,t}(\bar{x}_{t})\right]         \notag\\
&&~~~+2n{\color{black}m}\theta^{2}RLCG\alpha_{t}\sum_{l=1}^{t}\lambda^{t-l}\alpha_{l-1}  
+2R\Vert x_{t+1}^{\star}- x_{t}^{\star}\Vert
\notag\\
&&~~~+2n{\color{black}m}^{2}\theta^{2}CG^{2}\alpha_{t}\sum_{l=1}^{t+1}\lambda^{t+1-l}\alpha_{l-1}
+\frac{1}{2}{\color{black}m}^{2}\theta^{2}G^{2}\alpha_{t}^{2}  \notag\\
&&~~~+2\left(RL+{\color{black}m}G\lambda\right)\theta C\alpha_{t}\lambda^{t}\sum_{j=1}^{n}\Vert x_{j,0}\Vert 
\end{eqnarray}
Note that
\begin{eqnarray}\label{Expsssum}
-\sum_{t=0}^{T}\frac{\mathbb{E}[\Delta D_{t}]}{\alpha_{t}}
&\hspace*{-0.5em}=\hspace*{-0.5em}&
\sum_{t=0}^{T}\frac{\mathbb{E}[ D_{t}]-\mathbb{E}[ D_{t+1}]}{\alpha_{t}}   \notag\\
&\hspace*{-0.5em}\leq\hspace*{-0.5em}&
\frac{\mathbb{E}[ D_{0}]}{\alpha_{0}}
+\sum_{t=1}^{T}\left(\frac{1}{\alpha_{t}}-\frac{1}{\alpha_{t-1}}\right)\mathbb{E}[D_{t}]  \notag\\
&\hspace*{-0.5em}\leq\hspace*{-0.5em}&
\frac{2}{\alpha_{T}}R^{2}
\end{eqnarray}
where the positivity of $D_{t}$, $t\geq 0$, has been utilized to obtain the first inequality, and the last inequality follows from Assumption \ref{ConvexandBoundedG}.
{\color{black}Then, dividing both sides of (\ref{Expvarisu}) by $\alpha_{t}$ and summing the new inequality over $t\in[T]$, we obtain 
%
\begin{eqnarray}\label{Theorem22-Add}
&&\hspace{-1.5em}\sum_{t=0}^{T}\sum_{i=1}^{n}\frac{\pi_{i,t+1}}{\phi_{i,t+1}}\mathbb{E}\left[f_{i,t}(\bar{x}_{t})-f_{i,t}({\color{black}\tilde{x}_{t}^{\star}})\right]   \notag\\
&&\hspace{-1.5em}\leq
2n{\color{black}m}\theta^{2}RLCG\sum_{t=0}^{T}\sum_{l=1}^{t}\lambda^{t-l}\alpha_{l-1}
+\frac{{\color{black}m}^{2}\theta^{2}G^{2}}{2}\sum_{t=0}^{T}\alpha_{t}+\frac{2R^{2}}{\alpha_{T}}  \notag\\
&&\hspace{-1.5em}~~+2n{\color{black}m}^{2}\theta^{2}CG^{2}\sum_{t=0}^{T}\sum_{l=1}^{t+1}\lambda^{t+1-l}\alpha_{l-1}  
+\mu\theta G(1+T)  \notag\\
&&\hspace{-1.5em}~~+2\left(RL+{\color{black}m}G\lambda\right)\theta C\sum_{j=1}^{n}\Vert x_{j,0}\Vert\sum_{t=0}^{T}\lambda^{t} 
+\frac{2R}{\alpha_{T}}\mathcal{C}_{T}
\end{eqnarray}
%
Recalling the definition of $\tilde{x}_{t}^{\star}$, we have 
\begin{eqnarray}\label{regAADDD}
&&f_{i,t}(\bar{x}_{t})-f_{i,t}(x_{t}^{\star})\notag\\
&&=f_{i,t}(\bar{x}_{t})-f_{i,t}(\tilde{x}_{t}^{\star})+f_{i,t}(\tilde{x}_{t}^{\star})-f_{i,t}(x_{t}^{\star})\notag\\
&&\leq f_{i,t}(\bar{x}_{t})-f_{i,t}(\tilde{x}_{t}^{\star})+G\Vert  \tilde{x}_{t}^{\star}-x_{t}^{\star}\Vert\notag\\
&&\leq f_{i,t}(\bar{x}_{t})-f_{i,t}(\tilde{x}_{t}^{\star})+\xi GR
%
%
\end{eqnarray}
Finally, (\ref{Theorem22}) can be obtained by combining (\ref{Theorem22-Add}) and (\ref{regAADDD}). 
}

$\hfill \blacksquare$

\section*{Appendix D}

\textit{Proof of Theorem \ref{NetworkRegret}.}
Recalling (\ref{DynamicRegret}), we have 
\begin{eqnarray}\label{DynamicRegretTheo}
\mathcal{R}_{j}^{d}(T)
&\hspace*{-0.5em}=\hspace*{-0.5em}&
\sum_{t=0}^{T}\mathbb{E}\left[f_{t}(x_{j,t})\right]
-\sum_{t=0}^{T}f_{t}(x_{t}^{\star})   \notag\\
&\hspace*{-0.5em}=\hspace*{-0.5em}&
\sum_{t=0}^{T}\sum_{i=1}^{n}\mathbb{E}\left[f_{i,t}(x_{j,t})-f_{i,t}(\bar{x}_{t})\right]    \notag\\
&&+\sum_{t=0}^{T}\sum_{i=1}^{n}\mathbb{E}\left[f_{i,t}(\bar{x}_{t})-f_{i,t}(x^{\star}_{t})\right]
\end{eqnarray}
For the first term on the right side of (\ref{DynamicRegretTheo}), we have
\begin{eqnarray}\label{DynamicReSubxx}
f_{i,t}(x_{j,t})-f_{i,t}(\bar{x}_{t})
&\hspace*{-0.5em}\leq\hspace*{-0.5em}&
g_{i,t}^{\text{T}}(x_{j,t})(x_{j,t}-\bar{x}_{t})  \notag\\
&\hspace*{-0.5em}\leq\hspace*{-0.5em}&
G\Vert x_{j,t}-\bar{x}_{t}\Vert
\end{eqnarray}
where Assumption \ref{ConvexandBoundedG}(b) has been utilized to obtain the last inequality. 
Recalling (\ref{Difference1})-(\ref{wit}) and 
taking expectation on both sides of (\ref{DynamicReSubxx}), 
then 
the first term on the right  side of (\ref{DynamicRegretTheo}) can be bounded as follows
\begin{eqnarray}\label{DynamicReSub}
&&\sum_{t=0}^{T}\sum_{i=1}^{n}\mathbb{E}\left[f_{i,t}(x_{j,t})-f_{i,t}(\bar{x}_{t})\right] \notag\\
&&\leq
n^{2}{\color{black}m}\theta CG^{2}\sum_{t=0}^{T}\sum_{l=1}^{t}\lambda^{t-l}\alpha_{l-1}  \notag\\
&&~~~+nCG\sum_{i=1}^{n}\Vert x_{i,0}\Vert \sum_{t=0}^{T}\lambda^{t}.
\end{eqnarray}

To establish the desired dynamic regret bound, we now only need to bound the last term on the right side of (\ref{DynamicRegretTheo}).
Note that the term $\mathbb{E}[f_{i,t}(\bar{x}_{t})-f_{i,t}(x^{\star}_{t})]$ is scaled by $\frac{\pi_{i,t+1}}{\phi_{i,t+1}}$ in Theorem \ref{Lemma_4} and we cannot determine the sign of this term for particular $i\in\mathcal{V}$. 
Now, we calculate 
$-\frac{\pi_{i,t+1}}{\phi_{i,t+1}}(f_{i,t}(\bar{x}_{t})-f_{i,t}(x^{\star}_{t}))$ as follows 
\begin{eqnarray}\label{BasicRela1}
&&-\frac{\pi_{i,t+1}}{\phi_{i,t+1}}\left(f_{i,t}(\bar{x}_{t})-f_{i,t}(x^{\star}_{t})\right)  \notag\\
&&=
-\frac{\frac{1}{n}\phi_{i,t+1}+\pi_{i,t+1}-\frac{1}{n}\phi_{i,t+1}}{\phi_{i,t+1}}\left(f_{i,t}(\bar{x}_{t})-f_{i,t}(x^{\star}_{t})\right)   \notag\\
&&\leq
\frac{|\frac{1}{n}\phi_{i,t+1}-\pi_{i,t+1}|}{\phi_{i,t+1}}\big|f_{i,t}(\bar{x}_{t})-f_{i,t}(x^{\star}_{t})\big|   \notag\\
&&~~~-\frac{1}{n}\left(f_{i,t}(\bar{x}_{t})-f_{i,t}(x^{\star}_{t})\right)   \notag\\
&&\leq
\theta C\lambda^{T-t-1}G \Vert \bar{x}_{t}-x_{t}^{\star}\Vert
-\frac{1}{n}\left(f_{i,t}(\bar{x}_{t})-f_{i,t}(x^{\star}_{t})\right)
\end{eqnarray}
where the last inequality follows from Lemma \ref{BoundLe}-\ref{Rescale} and Assumption \ref{ConvexandBoundedG}(b).
Thus, we can further obtain
%
%
\begin{eqnarray}\label{BasicRela22}
&&\hspace{-1em}f_{i,t}(\bar{x}_{t})-f_{i,t}(x^{\star}_{t}) \notag\\
&&\hspace{-1em}\leq n\frac{\pi_{i,t+1}}{\phi_{i,t+1}}(f_{i,t}(\bar{x}_{t})-f_{i,t}(x^{\star}_{t}))  
+2n\theta CGR\lambda^{T-t-1}. 
\end{eqnarray}
Taking  expectation on both sides of (\ref{BasicRela22}) and combining Theorem \ref{Lemma_4}, we obtain   
\begin{eqnarray}\label{DynaRegretIm}
\mathcal{R}_{j}^{d}(T) 
&\hspace*{-0.8em}\leq\hspace*{-0.8em}&
n^{2}{\color{black}m}\theta CG(2\theta RL+G)\sum_{t=0}^{T}\sum_{l=1}^{t}\lambda^{t-l}\alpha_{l-1}  \notag\\
&&\hspace{-0.2em}+2n^{2}{\color{black}m}^{2}\theta^{2}CG^{2}\sum_{t=0}^{T}\sum_{l=1}^{t+1}\lambda^{t+1-l}\alpha_{l-1} 
+\frac{2nR^{2}}{\alpha_{T}}\notag\\
&&\hspace{-0.2em}+(2RL\theta+2{\color{black}m}\theta\lambda G+G)nC\sum_{i=1}^{n}\Vert x_{i,0}\Vert\sum_{t=0}^{T}\lambda^{t} \notag\\
&&\hspace{-0.2em}+2n^{2}\theta CGR\sum_{t=0}^{T}\lambda^{T-t-1}
+\frac{1}{2}n{\color{black}m}^{2}\theta^{2}G^{2}\sum_{t=0}^{T}\alpha_{t}   \notag\\
&&\hspace{-0.2em}+\frac{2nR}{\alpha_{T}}\mathcal{C}_{T}
+n\theta \mu G(T+1){\color{black}+n\theta \xi GR(T+1)}.
\end{eqnarray}

Now we turn to calculate the terms on the right side of (\ref{DynaRegretIm}), respectively.
For the first term, we have
\begin{eqnarray}\label{MB1}
\sum_{t=0}^{T}\sum_{l=1}^{t}\lambda^{t-l}\alpha_{l-1}
=
\sum_{t=0}^{T-1}\lambda^{t}\sum_{l=0}^{T-t-1}\alpha_{l}  \leq
\frac{1}{1-\lambda}\sum_{l=0}^{T-1}\alpha_{l}
\end{eqnarray}
Analogously, 
we obtain $\sum_{t=0}^{T}\sum_{l=1}^{t+1}\lambda^{t+1-l}\alpha_{l-1}\leq \frac{1}{1-\lambda}\sum_{l=0}^{T}\alpha_{l}$.
Note that 
$\sum_{t=0}^{T}\alpha_{t}\leq 1+\int_{0}^{T+1}(t+1)^{-\frac{1}{2}}dt \leq 2\sqrt{T+1}$. 
Then, (\ref{NetRegRes}) can be 
obtained,  
which completes the proof. 
$\hfill \blacksquare$
\bibliographystyle{IEEEtran}


\begin{IEEEbiography}
[{\includegraphics[width=1in,height=1.25in,clip,keepaspectratio]{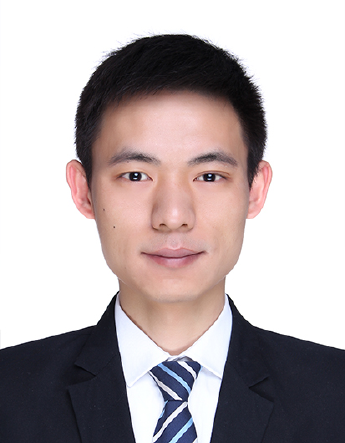}}]{Yongyang Xiong}
received the B.S. degree in information and computational science, the M.E. and Ph.D. degrees in control science and engineering from Harbin Institute of Technology, Harbin, China, in 2012, 2014, and 2020, respectively. From 2017 to 2018, he was a Joint Ph.D. Student with the School of Electrical and Electronic Engineering, Nanyang Technological University (NTU), Singapore. 
He is currently a Postdoctoral Fellow with the department of automation, Tsinghua University, Beijing, China. His current research interests include networked control system, distributed optimization and learning, multi-agent reinforcement learning, and their applications.

\end{IEEEbiography}


\begin{IEEEbiography}
[{\includegraphics[width=1in,height=1.25in,clip,keepaspectratio]{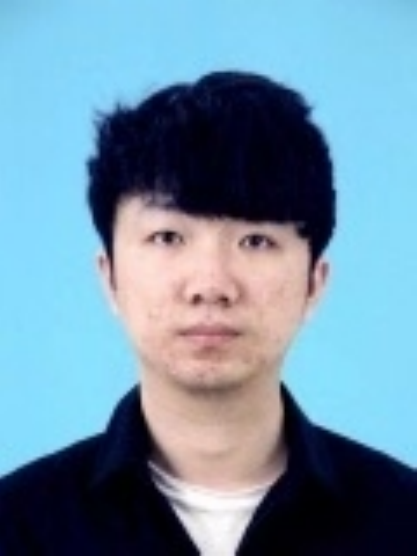}}]{Xiang Li}
received B.E. degree in Automation and M.E. degree in Control Science and Engineering from Harbin Institute of Technology, Harbin, China, in 2015 and 2018, respectively.
He is now pursuing his Ph.D. in the same major. His current research interests include image processing, pattern recognition, three-dimensional reconstruction and machine learning.
\end{IEEEbiography}


\begin{IEEEbiography}
[{\includegraphics[width=1in,height=1.25in,clip,keepaspectratio]{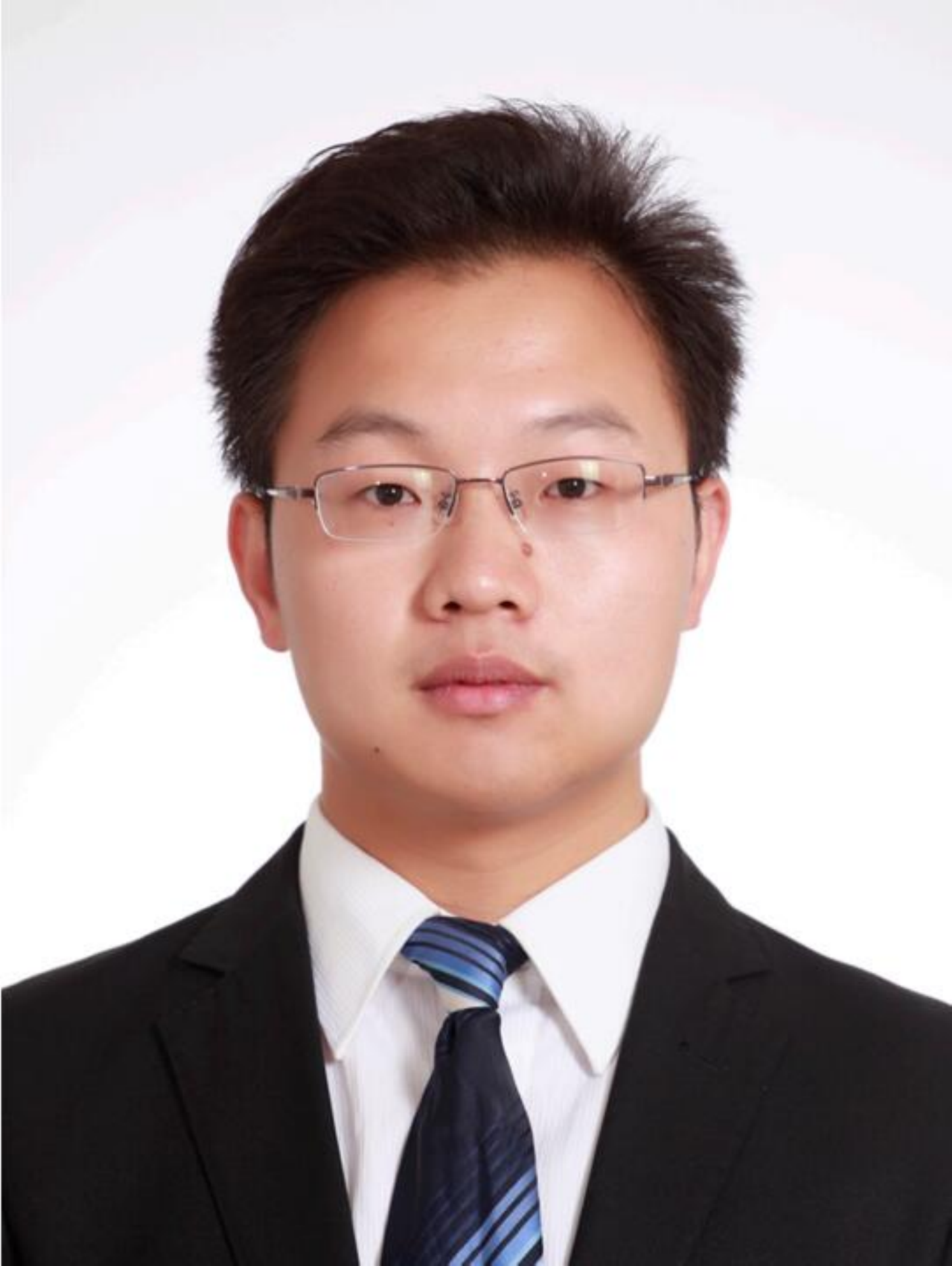}}]{Keyou You}
received the B.S. degree in statistical science from Sun Yat-sen University, Guangzhou, China, in 2007 and the Ph.D. degree in electrical and electronic engineering from Nanyang Technological University (NTU), Singapore, in 2012. After briefly working as a Research Fellow at NTU, he joined Tsinghua University, Beijing, China in 2012 where he is currently an Associate Professor with the Department of Automation. He held visiting positions at Politecnico di Torino, Turin, Italy, the Hong Kong University of Science and Technology, Hong Kong, and the University of Melbourne, Parkville, VIC, Australia. His current research interests include networked control systems, distributed algorithms, and their applications.

Dr. You was a recipient of the Guan Zhaozhi Award at the 29th Chinese Control Conference in 2010, the CSC-IBM China Faculty Award in 2014, and the National Science Fund for Excellent Young Scholars in 2017. He was nominated for the National 1000-Youth Talent Program of China in 2014.
\end{IEEEbiography}


\begin{IEEEbiography}
[{\includegraphics[width=1in,height=1.25in,clip,keepaspectratio]{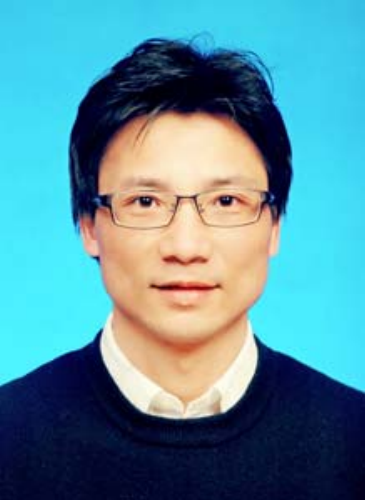}}]{Ligang Wu}(M'10-SM'12)
received the B.S. degree in Automation from Harbin University of Science and Technology, China in 2001; the M.E. degree in Navigation Guidance and Control from Harbin Institute of Technology, China in 2003; the PhD degree in Control Theory and Control Engineering from Harbin Institute of Technology, China in 2006. From January 2006 to April 2007, he was a Research Associate in the Department of Mechanical Engineering, The University of Hong Kong, Hong Kong. From September 2007 to June 2008, he was a Senior Research Associate in the Department of Mathematics, City University of Hong Kong, Hong Kong. From December 2012 to December 2013, he was a Research Associate in the Department of Electrical and Electronic Engineering, Imperial College London, London, UK. In 2008, he joined the Harbin Institute of Technology, China, as an Associate Professor, and was then promoted to a Full Professor in 2012. Dr. Wu was the winner of the National Science Fund for Distinguished Young Scholars in 2015, and received China Young Five Four Medal in 2016. He was named as the Distinguished Professor of Yangtze River Scholar in 2017, and was named as the Highly Cited Researcher in 2015, 2016 and 2017.

Dr. Wu currently serves as an Associate Editor for a number of journals, including \textsc{IEEE Transactions on Automatic Control}, \textsc{IEEE/ASME Transactions on Mechatronics}, \textit{Information Sciences}, \textit{Signal Processing}, and \textit{IET Control Theory and Applications}. He is also an Associate Editor for the Conference Editorial Board, IEEE Control Systems Society. Dr. Wu has published 6 research monographs and more than 150 research papers in international referred journals. His current research interests include switched systems, stochastic systems, computational and intelligent systems, multidimensional systems, sliding mode control, and flight control.
\end{IEEEbiography}

\balance
\end{document}